\documentclass[11pt]{article}
\usepackage{amsmath,amssymb,bbm}
\usepackage[dvips]{graphics}
\usepackage{graphicx}
\usepackage{psfrag}
\usepackage{graphics}
\usepackage{verbatim}
\RequirePackage[OT1]{fontenc}
\RequirePackage[latin1]{inputenc}

\title{High order expansion of matrix models and enumeration of maps.}

\author{
Edouard Maurel-Segala\thanks{Ecole Normale Sup\'erieure de Lyon,
Unit\'e de Math\'ematiques pures et appliqu\'ees,
UMR 5669,
46 All\'ee d'Italie, 
69364 Lyon Cedex 07, France. E-mail: emaurel@umpa.ens-lyon.fr}}

\newtheorem{prop}{Property}[section]
\newtheorem{theo}[prop]{Theorem}
\newtheorem{defi}[prop]{Definition}
\newtheorem{lem}[prop]{Lemma}

\newenvironment{dem}{\textbf{Proof.}\par}
{\begin{flushright}$\Box$\end{flushright}}

\def\e{\varepsilon}
\newcommand{\RR}{\mathbb{R}}
\newcommand{\CC}{\mathbb{C}}
\newcommand{\NN}{\mathbb{N}}

\newcommand{\HNC}{\mathcal{H}_N(\CC)}

\newcommand{\ppq}{\leqslant}
\newcommand{\pgq}{\geqslant}
\def\bt{{\bf t}}
\def\prob{\mathbb{P}}

\def\part{\partial}
\def\cycl{D}
\def\I{{\mathcal I}}
\def\oo{\overline}
\def\bA{{\bf A}}
\def\bB{{\bf B}}
\def\C{\mathcal C}
\def\MM{\mathcal M}
\def\Da{\mathcal D}

\def\a{\alpha}
\def\b{\beta}

\def\e{\varepsilon}

\def\l{\lambda}

\def\bk{{\bf k}}
\def\bj{{\bf j}}

\def\tr{\mathop{\rm tr}\nolimits}
\def\cxm{\CC\langle X_1,\cdots, X_m\rangle}

\def\bX{{\bf X}}

\def\hell{\hat{\ell}}

\def\munb{\oo{\mu}^N}
\def\nun{\nu^N}

\def\munv{\mu_V^N}
\def\munvt{\mu_{V_{\bt}}^N}

\def\mun{{\hat\mu^N}}
\def\sig2{{\sigma^2}}
\def\sur{\genfrac{}{}{0pt}{}}

\begin{document}
\maketitle
{\bf Abstract}
Perturbation of the {\bf GUE} are known in physics to be related to
enumeration of graphs on surfaces. Following \cite{GMa} and \cite{GMa2},
we investigate this idea and show that for a small convex perturbation,
we can perform a genus expansion: the moments
of the empirical measure can be developed into a series
whose $g$-th term is a generating function of graphs on a surface of genus
$g$.

\section{Introduction}
Wick's calculus allows to easily compute any moments of Gaussian variables
and gives them a combinatorial interpretation since the $p$-th moment of a
Gaussian
can be seen as the number of partitions in pairs of $[|1,p|]$.
This fact can be used to find moments of the $\bf{GUE}$,
the Gaussian
unitary model.
Let $\mu^N$ be the law on $\HNC^m$ the set of $m$-tuple 
$A_1,\cdots, A_m$ of $N\times N$ hermitian matrices
such that $\Re e A_i(kl),k<l$,
$\Im m A_i(kl),k<l$, $2^{-\frac{1}{2}}A_i(kk)$
is a family of independent real Gaussian of variance $(2N)^{-1}$ or more
directly
$$\mu^N(d\bA)=\frac{1}{Z^N}e^{{-\frac{N}{2}\tr(\sum_{i=1}^mA_i^2)}
}d^N\bA$$
with $d^N\bA$ the Lebesgue measure on $\HNC^m=(\RR^{N^2})^m$ and $Z^N$ 
a constant of normalization.

For a edge-colored graph on an orientated surface we say that a vertex is of
type 
$q=X_{i_1}\cdots X_{i_p}$ for a monomial $q$ if this vertex is of valence
$p$
and when we look at the half-edges going out of it, starting from a
distinguished one and going
in the clockwise order the first half-edge is of color $i_1$, the second
of color $i_2$,\dots, the $p$-th
of color $i_p$. A graph on a surface
is a map if it is connected and its faces are
homeomorphic to discs
(see section 3 for a precise definition of these notions).
Then, a simple computation (a simple generalization of \cite{zvon})
using Wick's calculus shows
that
for all non commutative monomials, $X_{i_1}\cdots X_{i_p}$,
$$\mu^N(\frac{1}{N}\tr(A_{i_1}\cdots A_{i_p}))=
\sum_{g\in \NN}\frac{1}{N^{2g}}\MM_g(X_{i_1}\cdots X_{i_p})$$
where $\MM_g(X_{i_1}\cdots X_{i_p})$ is the number up to isomorphism
of maps with colored edges
on a surface
of genus $g$ with one vertex of type $X_{i_1}\cdots X_{i_p}$.
Besides, one can use
Euler's formula to show that the sum in the right hand side is always finite.

Is such an interpretation can be generalized beyond the Gaussian case?
Take a potential
$V(X_1,\cdots,X_m)=\sum_i t_iq_i$.
with complex $t_1,\cdots,t_n$ and non-commutative monomials $q_i$.
We are interested in the following perturbation of the 
$\bf{GUE}$
\begin{equation}\label{guepert}\munv(dA_1,\cdots,dA_m)=\frac{1}{Z^N_V}
e^{-N\tr(V(A_1,\cdots, A_m))}
d\mu^N(A_1,\cdots,A_m)
\end{equation}
where $Z^N_V$ is the normalizing constant so that $\munv$ is
a probability measure. The derivatives of the moments of this model
in $\bt=0$ are exactly moments of the $\bf{GUE}$ and thus can be computed
using Wick's calculus and the limit can be formally expressed as
a generating function of graphs. For example, for $V=tX^4$, we can
obtain
the following expansion for the free energy
\begin{equation}\label{phi4}\frac{1}{N^2}\log Z^N_{tX^4}=\sum_{g\in\NN}
\frac{1}{N^{2g}}\sum_{k\in\NN}
\frac{(-t)^{k}}{k!}\C_g^{k}
\end{equation}
with $\C_g^{k}$ the number up to isomorphism of connected graphs on a
surface of genus
$g$ with $k$ vertex of valence $4$ and such that faces are homeomorphic
to discs (the so-called maps).
Note that we have to be careful since the right hand side of (\ref{phi4}) is
divergent for $t\neq 0$. Thus, this equality is purely formal but
we will be able to give it a precise mathematical meaning
(at least, Wick's calculus show that the derivatives of both
sides are equal in $t=0$).

Such genus expansion has been used for long in physics
and one may wonder for which potentials we can prove it.
In the one matrix case, the problem has been solved in \cite{EML}
using Riemann-Hilbert methods. This case has also been studied a lot in
physics, see for example the review \cite{FGZ}. The multi-matrix case
is much complicated since the technique of orthogonal polynomials
can not be generalized. Although, there is a large literature in physics
(see \cite{EKK}, \cite{BEH2} for example) on some
specific cases such
as the Ising model on random graphs i.e. a potential of the form
$V=V_1(A)+V_2(B)+cAB$. In mathematics, this problem has also been studied
with a completely different approach, namely large deviation in \cite{GCMP},
which gives the first asymptotic of models
closed to the Ising model.

In this paper we will try to avoid the use of orthogonal 
polynomials and our main tool will be the so-called Schynger-Dyson's 
equation. Besides, we will concentrate on the combinatorial
interpretation of the limit.
Our study follows the two papers \cite{GMa} and \cite{GMa2} which
respectively studied the first order asymptotic and a central limit theorem
for these models.

We will always assume
\begin{enumerate}
\item
the perturbation is small, we will restrict ourselves to small
coefficients $t_i$ in $V$. Note that we can not get rid of this condition
as the generating functions of combinatorial objects that appear have
arbitrary small radius of convergence.
\item
the potential $V+\frac{1}{2}\sum_iX_i^2$ is ``uniformly'' convex: there exists
$c>0$
such that for all $N$ in $\NN$,
$$\varphi^N_V:
\begin{array}{ccc}
\HNC^m&\longrightarrow&\CC\\
(X_1,\cdots,X_m)&\longrightarrow&\tr(
V(X_1,\cdots,X_m)+\frac{1-c}{2} \sum_{i=1}^m X_i^2)
\end{array}
$$
is a real and convex function. If $V$ satisfy this condition, we will
say that $V$ is $c$-convex.
\end{enumerate} 
Thus, for $V_{\bt}=\sum_{i=1}^nt_iq_i$ with $\bt=(t_1,\cdots,t_n)$ complex
numbers and $q_i$ non-commutative monomials we define
$$B_{\eta,c}=\{\bt\in\CC^n|
|\bt|=\max_i |t_i|\ppq \eta,\mbox{$V$ is $c$-convex}\}.$$
An example of $c$-convex function is 
$$V(\bX)=\sum_i  P_i(\sum_{j}\alpha_{ij}X_j)
+\sum_{k\ell}\beta_{k\ell}X_kX_\ell$$
with real and convex polynomials $P_i$, real $\alpha_{ij},\beta_{kl}$ and
for all $l$, $\sum|\beta_{kl}|<1-c$.

The main result of this paper is
\begin{theo}\label{theofree}
Let $V_{\bt}=\sum_{i=1}^nt_iq_i$, and $c>0$, for all $g\in \NN$, there exists
$\eta_g>0$
such that for all $\bt$ in
$B_{\eta_g,c}$, the free energy has the following expansion
$$F^N_{V_{\bt}}:=\frac{1}{N^2}\log Z^N_{V_{\bt}}=F^0(\bt)+\frac{1}{N^2}
F^1(\bt)+\cdots+\frac{1}{N^{2g}}F^g(\bt)+o(\frac{1}{N^{2g}})$$
with $F^g$ the generating function for maps of genus $g$ associated with $V$:
$$F^g(\bt)=
\sum_{\bk\in\NN^n}
\frac{(-\bt)^{\bk}}{\bk!}\C_g^{\bk}(P)$$
where $\bk!=\prod_i k_i!$, $(-\bt)^{\bk}=\prod_i(-t_i)^{k_i}$
and $\C_g^{\bk}$ is the number of maps on a surface
of
genus $g$ with
$k_i$ vertices of type $q_i$. 
\end{theo}
In addition to the free energy, we are interested in more general observables.
For example we would like to know the asymptotic of the non-commutative
moments of our measure
$E_{\mu^N_{V_{\bt}}}[\frac{1}{N}\tr(P)]$
for a non-commutative polynomial $P$. Such moments
appear as derivatives of the free energy since
$$-E_{\mu^N_{V_{\bt}}}[\frac{1}{N}\tr(P)]
=\left.\frac{\partial}{\partial u}\right|_{u=0}F^N_{V_{\bt}+uP}.$$
\begin{theo}\label{maintheo}
With the same hypothesis than in the previous theorem, for all
 for all $g\in \NN$, there exists $\eta>0$, 
such that for all $\bt$ in
$B_{\eta,c}$, for all monomials $P$
$$E_{\mu^N_{V_{\bt}}}[\frac{1}{N}\tr(P)]
=\C^0_{\bt}(P)+\cdots
+\frac{1}{N^{2g}}\C^g_{\bt}(P)+o(\frac{1}{N^{2g}})$$
with $\C_g$
 the generating function maps of genus $g$ with some fixed vertices:
$$\C^g_{\bt}(P)=
\sum_{\bk\in\NN^n}
\frac{(-\bt)^{\bk}}{\bk!}
\C_g^{\bk}(P)$$
where $\C_g^{\bk}(P)$ is the number of maps on a surface
of
genus $g$ with
$k_i$ vertices of type $q_i$ and one of type $P$.
\end{theo}

In fact, we will be able to find the asymptotic of much more observable
such as the higher derivatives of the free energy. Indeed, we show
that we can derivate term by term the expansion of Theorem \ref{theofree}.
Let us introduce for $\bj=(j_1,\cdots,j_n)\in \NN^n$,
the operator of derivation
$$\Da_\bj=\frac{\partial^{\sum_i j_i}}{\partial t_{1}^{j_1}\cdots
\partial t_{n}^{j_n}}.$$
\begin{theo}\label{theofree2}
With the same hypothesis than in the previous theorem, for all
$\bj=(j_1,\cdots,j_n)\in \NN^n$, for all $g\in \NN$, there exists $\eta>0$
such that for all $\bt$ in
$B_{\eta,c}$,
$$\Da_\bj F^N_{V_{\bt}}
=\Da_\bj F^0(\bt)+\cdots
+\frac{1}{N^{2g}}\Da_\bj F^g(\bt)+o(\frac{1}{N^{2g}}).$$
Besides, $\Da_\bj F^g$
is the generating function maps of genus $g$ with some fixed vertices:
$$\Da_\bj F^g(\bt)=
\sum_{\bk\in\NN^n}
\frac{(-\bt)^{\bk}}{\bk!}
\C_g^{\bk+\bj}.$$
\end{theo}

In the next section, we will define some useful notation when one has to 
deal with
non-commutative probability theory and we will recall the main result
of \cite{GMa}. Next, we will begin to look for recursive relations between
the asymptotic of the non-commutative moments of our model. This will lead
us  to study some combinatorial objects in section 4 whose generating
function satisfy these relations.
In the two next sections, we will prove the equality of
these moments and these enumerating
functions before
proving our main results. Finally the last section will be devoted
to the proof of Theorem \ref{theofree2}.
 
\section{Notations and reminder}
We denote by $\cxm$  the set of complex  polynomials
on the non-commutative
unknown $X_1$,\dots,$X_m$ i.e. the complex linear combination of monomials 
which are simply the set of finite words on $X_1,\ldots,X_m$. Monomials
must be thought as non-commutative moments.
 Let $*$ denotes the linear involution on $\cxm$
such that for all complex $z$ and all monomials
$$(zX_{i_1}\dots X_{i_p})^*=\oo{z}X_{i_p}\dots X_{i_1}.$$
We will say that a polynomial $P$ is self-adjoint if $P=P^*$.
We will denote $\cxm^*$ the dual of $\cxm$. An element of $\cxm^*$
has a compact support if there exists $R>0$
such that for all monomial
$|\tau(X_{i_1}\dots X_{i_p})|\ppq R^p$. By analogy with the one variable case,
the infimum of the $R$'s which satisfy this inequality for all monomials will
be called the radius of the support of $\tau$.

For a polynomial $P$ and a monomial $q$, we define $\lambda_q(P)$ as 
the coefficient of $q$ in the decomposition of $P$.
For $M>0$, we define the norm $\|.\|_M$ on polynomials:
$$||P||_M=\sum_{l\in \NN}\sum_{\sur{q\textrm{ monomial}}
{deg q=l }}|\lambda_q(P)|M^l.$$
This norm $\|.\|_M$ is an algebra norm, i.e. for all polynomials
$P$, $Q$, $$\|PQ\|_M\ppq \|P\|_M\|Q\|_M.$$
Note that an element $\tau$ of $\cxm^*$ has support of radius less than $R$
if and only if for all polynomials $P$,
$$|\tau(P)|\ppq \|P\|_R.$$
We extend $\|.\|_M$ on $\cxm\otimes\cxm$ by defining this norm on the decomposition
in monomials:
$$\|\sum\l_{q_1,q_2}q_1\otimes q_2\|_M=\sum|\l_{q_1,q_2}|M^{\deg q_1+\deg q_2}$$
with this definition for all polynomials $P,Q$, $\|P\otimes Q\|_M=\|P\|_M\|Q\|_M$.

We define the 
non-commutative derivatives $\part_i$ from $\cxm$ to $\cxm^{\otimes 2}$
for $1\ppq i\ppq m$ by the Leibniz rule
$$\part_i PQ=\part_iP\times(1\otimes Q)+(P\otimes 1)\times \part_iQ$$
and 
$\part_i X_j={\mathbbm 1}_{i=j}1\otimes 1.$
For a monomial $P$, we will often use the convenient expression
$$\part_iP=\sum_{P=RX_iS}R\otimes S$$
where the sum runs over all possible monomials $R,S$
so that  $P$ decomposes into $RX_iS$.
We also define another operator of derivation on polynomials, the cyclic
derivative $\cycl_i$ which
is linear and such that for all monomials:
$$\cycl_iP=\sum_{P=RX_iS}  SR.$$
Alternatively, $\cycl$ can be defined as $m\circ \part$ where
$m(A\otimes B)=BA$.
We will see that these two operators appear naturally when we derivate
products of matrices and they both possesses a nice combinatorial
interpretation. An important fact we will use later is that for all $M'>M$, both
$\part_i$ and $\cycl_i$ are continuous from $(\cxm,\|.\|_{M'})$ to $(\cxm,\|.\|_{M})$.
For example for a monomial $q$,
$$\frac{\|\cycl_i q\|_M}{\|q\|_{M'}}=\deg q \frac{M^{\deg q-1}}{{M'}^{\deg q}}= M^{-1}\deg q(\frac{M}{M'})^{\deg q}$$
which is bounded.

Our main object of study is the law $\munvt$ on $\HNC^m$
$$
\munvt(dA_1,\cdots,dA_m)=\frac{1}{Z^N_{V_{\bt}}}
e^{-N\tr(V_{\bt}(A_1,\cdots, A_m))}
d\mu^N(A_1,\cdots,A_m)
$$
and we are particulary interested by the behavior of the
random variable
$$\mun:
\begin{array}{ccc}
\cxm&\longrightarrow&\CC \\
P&\longrightarrow&\frac{1}{N}\tr(P(A_1,\cdots,A_m))
\end{array}
$$
and its mean:
$$\munb_{\bt}:
\begin{array}{ccc}
\cxm&\longrightarrow&\CC \\
P&\longrightarrow&E_{\munvt}[\frac{1}{N}\tr(P(A_1,\cdots,A_m))]
\end{array}
$$

We can now state precisely the main result of \cite{GMa}, we will overuse
it in the next sections.
For any $c$, there exists $\eta>0$ such that for all $\bt\in B_{\eta,c}$,
for all polynomials $P$,
$\mun(P)$
goes when $N$ goes to $+\infty$, almost surely and in expectation towards
$\mu_{\bt}(P)$ with
$\mu_{\bt}$ a
solution of the Schwinger-Dyson equation
$$\mu_{\bt}\otimes \mu_{\bt}(\part_iP)=
\mu_{\bt}((X_i+\cycl_i V_{\bt})P)\quad\forall P\in\cxm, \,
\forall i\in\{1,\cdots,m\}.$$
Besides, this soluion has a bounded support, there exists $R$ such that
for all monomial $X_{i_1}\cdots X_{i_p}$
\begin{equation}\label{support}
|\mu_{\bt}(X_{i_1}\cdots X_{i_p})|\ppq R^p.
\end{equation}
and $\mu_{\bt}$ can be seen as a generating function
of planar maps, for all polynomials $P$,
$$\mu_{\bt}(P)=\sum_{k_1,\cdots,k_n\in\NN}
\prod_i\frac{(-t_i)^{k_i}}{k_i!}\MM^{k_1,\cdots,k_n}(P)$$
where $\MM^{k_1,\cdots,k_n}(P)$ is the number of planar connected graphs with
$k_i$ vertices of type $q_i$ and one of type $P$. For the rest of the paper
we will work in this domain $B_{\eta,c}$ where the convergence holds.
In order to shorten a little the notations the subscript $\bt$ will be most
of the time
implicit, for example we will often write $\mu$ instead of $\mu_{\bt}$,
$V$ instead of $V_{\bt}$, $\munb$ instead of $\munb_{\bt}$ \dots
\section{First order observable}
The starting point is a relation already used
in \cite{GMa} for the
matrix model when $N$ is fixed: for all polynomial $P$, for all $i$,
$$E[\mun((X_i+\cycl_i V) P)]=E[(\mun\otimes\mun)(\part_i P)].$$
We will give the proof of a generalization of this equality later.
Using this equality and some concentration inequalities we were able to prove
that for $\bt$ in $B_{\eta,c}$
for a well chosen $\eta$, for all polynomial $P$,
$E[\mun(P)]$ was converging towards $\mu(P)$ with $\mu$
the unique solution of the Swinger-Dyson's equation {\bf SD[V]}:
\begin{equation}\label{finitesd}
\mu((X_i+\cycl_i V) P))=(\mu\otimes\mu)(\part_i P).
\end{equation}
In order to find the next asymptotic, we study the difference between the
equation
for finite $N$ and the limit equation, if $\nun=N^2(\munb-\mu)$, we obtain by
substracting
the two equations:
\begin{equation}\label{premierecorrec}
\nun((X_i+\cycl_i V) P)-(I\otimes\mu+\mu\otimes I)\part_i P)=N^2
E[((\mun-\mu)\otimes(\mun-\mu))(\part_i P)]
\end{equation}
Here, an important operator shows up in the left hand side.
Let $P=X_{i_1}\cdots X_{i_p}$ be a monomial and define the following
operators:
$$\Xi_1 P=\frac{1}{p}\sum_i \cycl_i V \cycl_i P$$
$$\Xi_2 P=\frac{1}{p}\sum_i (I\otimes\mu+\mu\otimes I)\part_i \cycl_i P.$$
We extend them by linearity on $\cxm$ and we define $\Xi_0=I-\Xi_2$ and
$\Xi=\Xi_0+\Xi_1$. 
These operators were introduced in \cite{GMa2} to obtain a central limit
theorem for the matrix model. We also define the operator of division by the
degree i.e. the
linear operator $P\to\oo{P }$
such that for all monomial $P=X_{i_1}\cdots X_{i_p}$,
$\oo{P}=\frac{1}{p}P$ and $\oo{1}=0$.
These operators allow us to state the relation for the first correction
(\ref{premierecorrec}) in a simpler form:
\begin{equation}\label{equ1}
\nun(\Xi P)=N^2E[((\mun-\mu)\otimes(\mun-\mu))(\part_i \cycl_i \oo{P})].
\end{equation}
Then, the strategy is simple, we only have to understand the asymptotic of 
$N^2E[((\mun-\mu)\otimes(\mun-\mu))(R\otimes S)]$ and then "invert" $\Xi$.
The first order asymptotic of 
$N^2E[((\mun-\mu)\otimes(\mun-\mu))(R\otimes S)]$
is easy to compute using \cite{GMa2} as it was shown that $N(\mun-\mu)$
converges in law towards
a Gaussian when $N$ goes to infinity. The main issue is that when we try
to investigate the next asymptotic,
terms of type $N^3E[(\mun-\mu)^{\otimes 3}(R\otimes S\otimes T)]$ will
appear and at their turn they will create terms
of greater complexity. That's why we are interested more
generally in all the
$N^{\ell}E[(\mun-\mu)^{\otimes \ell}(P_1\otimes \cdots\otimes P_\ell)]$'s and
we will eventually
find  their full asymptotic.
First remark that according to \cite{GMa2}, for all $P$, $N(\mun(P)-\mu(P))$ 
converges to a gaussian
variable and this convergence occurs in moments (see Corollary 4.8
in \cite{GMa2}).
Thus $N^{\ell}E[(\mun-\mu)^{\otimes \ell}(P_1\otimes \cdots\otimes P_\ell)]$
has a finite limit when $N$ goes to infinity and this limit is $0$ 
if $\ell$ is odd. But we need a more precise result which state that this
convergence is uniform for all monomials
$P$ of reasonable degree.
\begin{lem}\label{estimate}
For all $\ell\in\NN^*$, $\a>0$ there exists $C,\eta,M_0>0$,
such that for all $\bt\in B_{\eta,c}$, $M>M_0$ and
all polynomials $P_1$,\dots,$P_\ell$ of degree less than $\a N^{\frac{1}{2}}$,
$$|E[N^{\ell}(\mun-\mu)^{\otimes \ell}
(P_1\otimes\cdots\otimes P_{\ell})]|\leq C\|P_1\|_M\cdots\|P_\ell\|_M.$$
\end{lem}
\begin{dem}
First using Hölder's inequality, write
$$E[N^{\ell}(\mun-\mu)^{\otimes \ell}
(P_1\otimes\cdots\otimes P_{\ell})]\leq\prod_{r=1}^\ell E[
|N(\mun-\mu)
(P_r)|^\ell]^{\frac{1}{\ell}}.$$
Thus we only have to prove the claim if the $P_r$ are equals.
Then we substract the mean
\begin{equation}\label{eqdec}
E[|N(\mun-\mu)(P)|^\ell]\leq 2^\ell\{E[|N(\mun(P)-\munb(P)|^\ell]
+|N(\munb(P)-\mu(P))|^\ell\}
\end{equation}
Proposition 3.1 in \cite{GMa2} state a rate of convergence of $\munb(P)$
to $\mu(P)$: there exists $C,M_0>0$ such that for $M>M_0$, for
all polynomials $P$
of degree less than $\e N^{\frac{2}{3}}$,
\begin{equation}\label{eqmoy}
|N(\munb(P)-\mu(P))|\ppq C\frac{\|P\|_M}{N}.
\end{equation}
Thus in inequality (\ref{eqdec}) above,  we only have to control the first term.
As we restrict ourself on the domain where the support of $\mu$ is
uniformly bounded by $R$, according to (\ref{support})
we get that for all polynomials of degree
less than $N^{\frac{1}{2}}$, 
$$
|\munb(P)|\ppq |\mu(P)|+C\frac{\|P\|_M}{N}\ppq C' \|P\|_M
$$
for $M$ larger than $M_0,R$.
We split the first term in two, let
$\mun_-=\mun{\mathbbm 1}_{\|\bA\|\leq M}$ and 
$\mun_+=\mun{\mathbbm 1}_{\|\bA\|> M}$ with $\|.\|$ the operator-norm.
By Jensen's inequality
$$E[|N(\mun_+(P)-E[\mun_+(P)]|^\ell]
\ppq (2N)^\ell\munb((PP^*)^\ell)^{\frac{1}{2}}
\prob(\|\bA\|> M)^{\frac{1}{2}}
$$
Now according to Lemma 2.2 in \cite{GMa2}, there exists $\alpha>0$ such
that if $M$ is sufficiently large,
$\prob(\|\bA\|> M)\ppq e^{-\alpha MN}$.
Thus
\begin{equation}\label{eqplus}
E[|N(\mun_+(P)-E[\mun_+(P)]|^\ell]\ppq C\|(PP^*)^\ell\|_M^{\frac{1}{2}}
\ppq C\|P\|_M^\ell
\end{equation}
Observe that $\mun(P)$ is lipschitz on $\{\|\bA\|\ppq M\}$ with lipschitz constant
\begin{align*}
&\sup_{\{\|\bA\|\ppq M\}}(\sum_{i,\a,\b}|\frac{\partial}
{\partial A_i(\a\b)}\frac{1}{N}\tr P(\bA)|^2)^{\frac{1}{2}}\\
&= \sup_{\{\|\bA\|\ppq M\}}(\sum_{i,\a,\b} |\frac{1}{N}(\cycl_iP(\bA))_{\a\b}|^2)^{\frac{1}{2}}\\
&= \sup_{\{\|\bA\|\ppq M\}}( \sum_i\frac{1}{N^2}\tr(\cycl_iP(\cycl_iP)^*))^{\frac{1}{2}}\\
&\ppq \frac{C}{\sqrt{N}}\sup_i\|\cycl_iP\|_M\ppq \frac{C'}{\sqrt{N}}\|P\|_{M'}
\end{align*}
where we choose $M'>M$ and use the continuity of the cyclic derivative.
We want to use a concentration inequality, thus we define an extension $\varphi$
of $\mun(P)$ on $\HNC^m$ which is equal to $\mun(P)$ on $\{\|\bA\|\ppq M\}$
and with the same lipschitz constant. For example, one can define:
$$\phi(\bA)=\sup_{\|\bB\|\ppq M}
(\frac{1}{N}\tr(P(\bB))-\frac{C'}{\sqrt{N}}\|P\|_{M'}\|\bA-\bB\|_2).$$
It is an easy exercise to check that $\phi$ is a
$\frac{C'}{\sqrt{N}}\|P\|_{M'}$-lipschitz function
which coincide with $\mun(P)$ on $\{\|\bA\|\ppq M\}$.
Now observe that with our hypothesis, for all $\bt$ in $B_{\eta,c}$, $\munv$ satisfy
a log-sobolev inequality with constant $Nc$. Thus, according to Herbst's argument (see \cite{ane})
we obtain a concentration's inequality, for all $\l$-lipschitz function $f$,
$$\munv(|f-\munv(f)|>\e)\ppq 2e^{-\frac{Nc\e^2}{\l^2}}.$$
Therefore,
$$\munv[|N(\mun_-(P)-E[\mun_-(P)]|>\e]\ppq 2e^{-\frac{C\e^2}{\|P\|_{M'}^2}}$$
and $$E[|N(\mun_-(P)-E[\mun_-(P)]|^\ell]\ppq C'\|P\|_{M'}^\ell.$$
Thus with (\ref{eqplus})
$$E[|N(\mun(P)-E[\mun(P)])|^\ell]\leq 2^\ell(C'+C)\|P\|_M^\ell$$
and we conclude with (\ref{eqdec}) and (\ref{eqmoy}). 
\end{dem}
We now try to find some relation between the
$N^{\ell}E[(\mun-\mu)^{\otimes \ell}(P_1\otimes \cdots\otimes P_\ell)]$'s
which generalize (\ref{equ1}). Remember that, for $\ell$ odd, those
quantities vanish when $N$
goes to infinity. Thus in order to obtain non-trivial limits, we have
to distinguish the 
normalisation according to the parity of $\ell$.

For $\ell\geq 1$ we define: 
$$\hat{\ell}=\left\{
\begin{array}{ll}
\ell+1&\mbox{ if $\ell$ is odd}\\
\ell&\mbox{ otherwise.}
\end{array}\right.
$$
Note that $\hat{\ell}$ is always an even integer.
We now define a function from $\sqcup_{\ell\in\NN} \cxm^{\otimes \ell}$
to $\CC$,
$$\nu^{N}(P_1\otimes\cdots\otimes P_\ell)=E_{\mu_V^N}[
N^{\hat{\ell}}(\mun-\mu)^{\otimes \ell}
(P_1\otimes\cdots\otimes P_\ell)].$$
On $\cxm^{\otimes \ell}$, this is a $\ell$-linear symmetric function
which is tracial in each $P_r$.
Our convention will be that for $\lambda$ in $\CC=\cxm^{\otimes 0}$,
 $\nu^N(\l)=\l$.
The relation that will appear as our main tool are the aim of the next
property. In a tensor product
$P_1\otimes \cdots \check{P_r} \cdots\otimes P_\ell$ denotes the tensor
product of  $P_1$,\dots,$P_{r-1}$,$P_{r+1}$,
\dots,$P_\ell$ i.e. the term $P_r$ is omitted.
\begin{prop}\label{sd1}
For all $\ell$, for all polynomial $P_1$,\dots,$P_\ell$, for all $N$,
if $\ell$ is even
\begin{align*}
\nu^N(\Xi P\otimes P_2\otimes \cdots\otimes P_\ell)
&=\sum_{i,r}\mu_{\bt}(\cycl_i\oo{P}\cycl_i P_r)\nu^N
( P_2\otimes \cdots \check{P_r} \cdots\otimes P_\ell)\\
&+\frac{1}{N^2}\sum_{i,r}\nu^N
(\cycl_i\oo{P}\cycl_i  P_r\otimes P_2\otimes \cdots \check{P_r} \cdots\otimes
P_\ell)\\
&+\frac{1}{N^2}\sum_i\nu^N
(\part_i\cycl_i\oo{P}\otimes P_2\otimes\cdots\otimes P_\ell)
\end{align*}
and if $\ell$ is odd
\begin{align*}
\nu^N(\Xi P\otimes P_2\otimes \cdots\otimes P_\ell)
&=\sum_{i,r}\mu_{\bt}(\cycl_i\oo{P}\cycl_i P_r)\nu^N
( P_2\otimes \cdots \check{P_r} \cdots\otimes P_\ell)\\
&+\sum_{i,r}\nu^N
(\cycl_i\oo{P}\cycl_i  P_r\otimes P_2\otimes \cdots \check{P_r} \cdots\otimes
P_\ell)\\
&+\sum_i\nu^N
(\part_i\cycl_i\oo{P}\otimes P_2\otimes\cdots\otimes P_\ell)
\end{align*}
\end{prop}
This property is the generalization of the equation (\ref{equ1}). 
One may wonder why we stress so much the difference between the odd and
the even case.
The point is to keep in mind which terms are of order $1$ and which are
negligible.
In view of this, the $\nu_1^N$ are convenient as they should all be of
order $1$ and thus
the previous equation will lead us to find their limit by induction.

\begin{dem}
To sum up the property in a shorter way, we have to prove that for all $\ell$
for all polynomials $P_1$,\dots,$P_\ell$ and for all $N$,
\begin{align*}
&N^\ell E[(\mun-\mu)^{\otimes \ell}(\Xi P_1\otimes P_2\otimes
\cdots\otimes P_\ell)]\\
&=\sum_{i,r}\mu_{\bt}(\cycl_i\oo{P_1}\cycl_i P_r)N^{\ell-2}
E[(\mun-\mu)^{\otimes \ell-2}
( P_2\otimes \cdots \check{P_r} \cdots\otimes P_\ell)\\
&+\sum_{i,r}N^{\ell-2}E[(\mun-\mu)^{\otimes \ell-1}
(\cycl_i\oo{P_1}.\cycl_i  P_r\otimes P_2\otimes \cdots \check{P_r}
\cdots\otimes P_\ell)]\\
&+\sum_iN^{\ell}E[(\mun-\mu)^{\otimes \ell+1}
(\part_i\cycl_i\oo{P_1}\otimes P_2\otimes\cdots\otimes P_\ell)].
\end{align*}
We will use the integration by part formula:
$$\int xf(x)e^{-x^2/2}dx=\int f'(x)e^{-x^2/2} dx.$$
We generalize this formula into
\begin{align*}
\int\tr(A_iP)f(\tr Q)d\mu^N&=\frac{1}{N}\sum_{\a,\b}\int
(\partial_{A_i(\a\b)}P_{\b\a})f(\tr Q)\\
&+P_{\b\a}
\partial_{A_i(\a\b)}
\tr Q f'(\tr Q)d\mu^N.
\end{align*}
Two useful computations show the importance of the non-commutative
derivatives and their links
with the derivation of polynomials of hermitian matrices: if $P$
is a monomial,
$$\sum_{\a\b}\partial_{A_i(\a\b)}P_{\b\a}=\sum_{\a\b}
\sum_{P=RX_iS}R_{\b\b}S_{\a\a}=\tr\otimes\tr(\part_i P)$$
and
$$\partial_{A_i(\a\b)}\tr P=\sum_{P=RX_iS,\gamma}R_{\gamma\b}
S_{\a\gamma}=(\cycl P)_{\a\b}.$$
Thus, for $P_1$,\dots,$P_\ell$ polynomials:
\begin{align*}
&N^\ell E[\mun(X_i P_1)(\mun-\mu)^{\otimes \ell-1}(P_2\otimes
\cdots\otimes P_\ell)]\\
&=N^{\ell}E[(\mun\otimes\mun)(\part_i P_1)(\mun-\mu)^{\otimes
\ell-1}(P_2\otimes \cdots\otimes P_\ell)]\\
&+\sum_r  N^{\ell-2}E[\mun(P_1\cycl_iP_r)(\mun-\mu)^{\otimes \ell-1}
(P_2\otimes \cdots \check{P_r}\cdots\otimes P_\ell)]\\
&- N^{\ell}E[\mun(P_1\cycl_iV)(\mun-\mu)^{\otimes \ell-1}(P_2\otimes
\cdots\otimes P_\ell)].
\end{align*}
Now remember that according to Schwynger Dyson's equation we have:
$$
\mu((X_i+\cycl_i V)P_1)
-\mu\otimes \mu(\part_iP_1)=0$$
Then, we substract the two equalities and use the identity 
$$\mun\otimes\mun-\mu\otimes\mu=(\mun-\mu)(I\otimes\mu+\mu\otimes I)+
(\mun-\mu)\otimes (\mun-\mu)$$
to obtain
\begin{align*}
&N^\ell E[(\mun-\mu)((X_i +\cycl_iV)P_1-(I\otimes\mu+\mu\otimes I)\part_i
P_1)\\
&\phantom{N^\ell E[}(\mun-\mu)^{\otimes \ell-1}(P_2\otimes \cdots\otimes
P_\ell)]\\
&=N^{\ell}E[(\mun-\mu)^{\otimes \ell+1}(\part_iP\otimes P_2\otimes \cdots
\otimes P_\ell)]\\
&+\sum_r  N^{\ell-2}E[(\mun-\mu)^{\otimes \ell-1}(P\cycl_iP_r\otimes P_2
\otimes \cdots \check{P_r}\cdots\otimes P_\ell)]\\
&+\sum_r  N^{\ell-2}\mu(P\cycl_iP_r)E[(\mun-\mu)^{\otimes \ell-2}\otimes 
P_2\otimes\cdots \check{P_r}\cdots\otimes P_\ell)].
\end{align*}
To get the result it is now sufficient to aplly this equality with
$P_1=\cycl_i {\bar P}$ and then to sum on $i$.
\end{dem}
This property gives us some precious hints on the limit $\nu$ of the $\nu^N$.
It should satisfy the "limit equation",
if $\ell$ is even
\begin{equation}\label{limiteven}
\nu(\Xi P\otimes P_2\otimes \cdots\otimes P_\ell)
=\sum_{i,r}\mu_{\bt}(\cycl_i\oo{P}\cycl_i P_r)\nu
( P_2\otimes \cdots \check{P_r} \cdots\otimes P_\ell)
\end{equation}
and if $\ell$ is odd
\begin{align}\label{limitodd}
\nu(\Xi P\otimes P_2\otimes \cdots\otimes P_\ell)
&=\sum_{i,r}\mu_{\bt}(\cycl_i\oo{P}\cycl_i P_r)\nu
( P_2\otimes \cdots \check{P_r} \cdots\otimes P_\ell)\\
&+\sum_{i,r}\nu
(\cycl_i\oo{P}\cycl_i  P_r\otimes P_2\otimes \cdots \check{P_r} \cdots\otimes
P_\ell)\nonumber\\
&+\sum_i\nu
(\part_i\cycl_i\oo{P}\otimes P_2\otimes\cdots\otimes P_\ell).\nonumber
\end{align}
Hopefully, we will be able to study the solutions $\nu$ of these equations.
In fact, following \cite{GMa} we would be able to prove that for $R,L>0$, 
there exists $\e>0$ such that for $|\bt|<\e$
there exists an unique $\nu:\sqcup_{\ell=0}^{2L}\cxm^{\otimes \ell}\to\CC$
linear on each
set $\cxm^{\otimes \ell}$ with support bounded by $R$
and which satisfy each of the previous equation for $\ell\leq 2L$.
But we will proceed in a different way. Looking at these equations, we will
try to 
recognize in them some relations between enumeration of combinatorial
objects. This is
the aim of the next section.
\section{Maps of high genus}
In this section, we describe the combinatorials
objects that appear in the asymptotic of our measure.
Remember that it was shown in \cite{GMa} that the first asymptotic can
be viewed as a generating function for the enumeration of planar 
maps with vertices of a given type.

First, we choose $m$ colors $\{1,\cdots,m\}$, one for each variable $X_i$.
A star must be thought as the neighbourhood of a vertex in a plane graph.
More precisely, it is a vertex with the half-edges coming out of it.
One of these
half-edge is distinguished and starting from it the other one are
clockwisely ordered.
Besides, each of these half-edges is colored.

We say that a star is of type $q$ for a monomial $q=X_{i_1}\cdots X_{i_p}$
if it has $p$ half-edges, the first half-edge
is distinguished and of color $i_1$ and then in the clockwise
order the second half-edge
is of color $i_2$, the third of color $i_3$, \dots, the
$p$-th of color $i_p$.
This gives a bijection between monomials and stars.

The combinatorial objects that will appear in the asymptotic of our
matrix model
are maps. A map is a connected graph on a compact orientated connected
surface such that
edges do not cross each other and faces are homeomorphic to discs.
We will consider edge-colored maps such that each vertex as a
distinguished edge going out of it so that we can associate a star and a well
defined type to
any vertex. 
The genus of the map is the genus of the surface.
We will count maps up to homeomorphism of the surface which preserves the
graph.

The typical way to construct a map is to put some stars $q_1,\cdots,q_p$
on a surface of genus $g$.
Then we consider all the half-edges that goes outside the stars
and glue them two by two while respecting the following constraints:
\begin{itemize}
\item
Two half-edges can only be glued if they are of the same color
\item
The edges created by gluing two half-edges mustn't cross any other edge.
\item
At the end of the process faces must be homeomorphic to discs.
\end{itemize} 

\begin{figure}[ht!]
\psfrag{pointint}{\large{?}}
\psfrag{x1}{$X_1$}
\psfrag{x2}{$X_2$}
\psfrag{fs}{first star}
\psfrag{ss}{second star}
\psfrag{de}{distinguished edges}
\begin{center}{\includegraphics{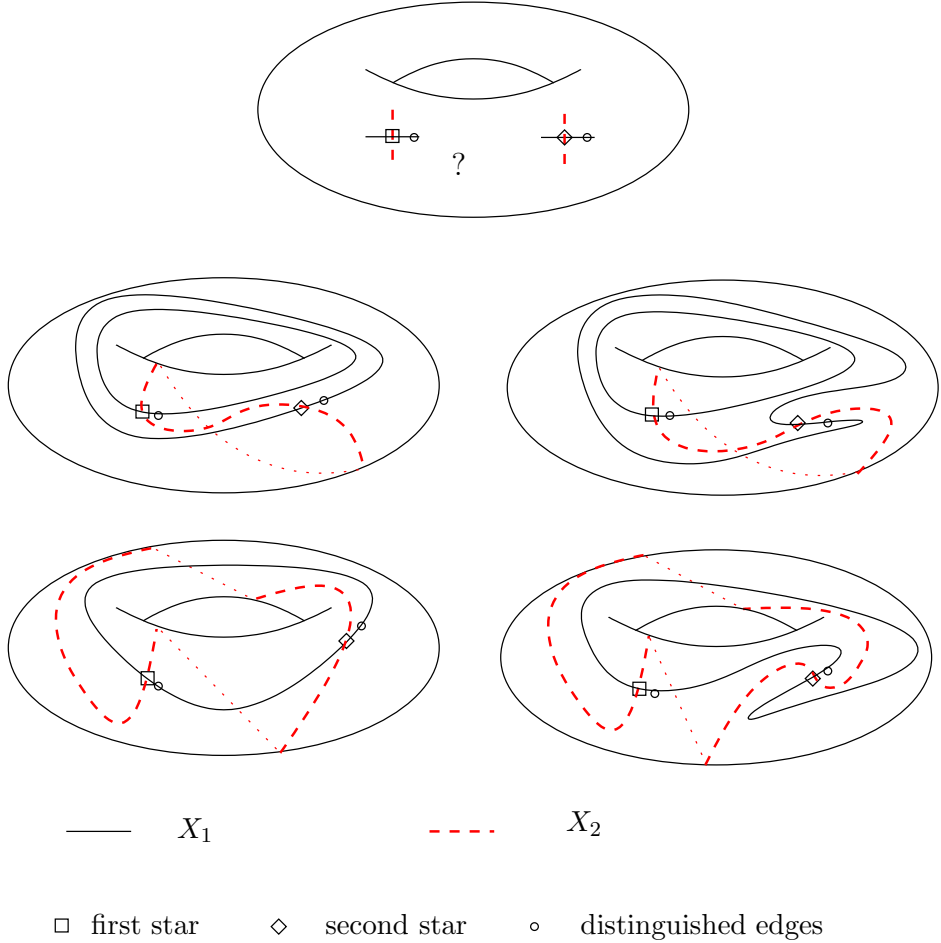}}
\end{center}
\caption{Maps of genus $1$ above two stars of type $X_1X_2X_1X_2$.}\label{tore}
\end{figure}

For example, one can ask how many maps of genus $1$, we can construct 
above two
stars of
type $X_1X_2X_1X_2$. The answer as shown in figure \ref{tore} is $4$. Note that 
as faces are
homeomorphic
to discs, it is sufficient to know which pairs of half-edges are glued
together to build the map.
\begin{defi}
For $\ell$ in $\NN$, $P$ a monomial
and integers $\bk=(k_1,\cdots,k_n)$, let
$\MM^{\bk}(P)$ be the number of 
maps of genus $0$ (or planar maps) with for all $i$, $k_i$ stars of
type $q_i$ and one of type $P$
where $q_1,\cdots,q_n$ are the monomials 
which appear in the potential $V$.
\end{defi}
We could define the same kind of quantities with the condition
of being of genus $g$ for $g>1$ but then we won't be able to find any
closed relation of induction between these quantities. In order to get relation
induction on enumeration of maps we follow an idea of Tutte (see \cite{Tu}).
We try to
decompose
a map in smaller ones by contracting one edge (Note that Tutte used to work
on the dual of the graph we are considering, thus his operation
is a little different).

Imagine a map of genus $1$ with a root of type $P=XRXS$ and that the two
half-edges corresponding
to the $X$ are glued together. Imagine also that the loop
resulting from this operation
is not retractable on the surface. How does the contraction of this edge
decomposes the map? Now $R$ and $S$
are separated by that loop, we will have to remind these two monomials.
That's why we will introduce maps
above a root of type $R\otimes S$. Besides $R$ and $S$ must be linked
together, otherwise there would be a face
(touched by the loop) which is not a disc, something to avoid for a map.

\begin{figure}[ht!]
\psfrag{x1}{$X_1$}
\psfrag{x2}{$X_2$}
\psfrag{de}{distinguished edge}
\begin{center}{\includegraphics{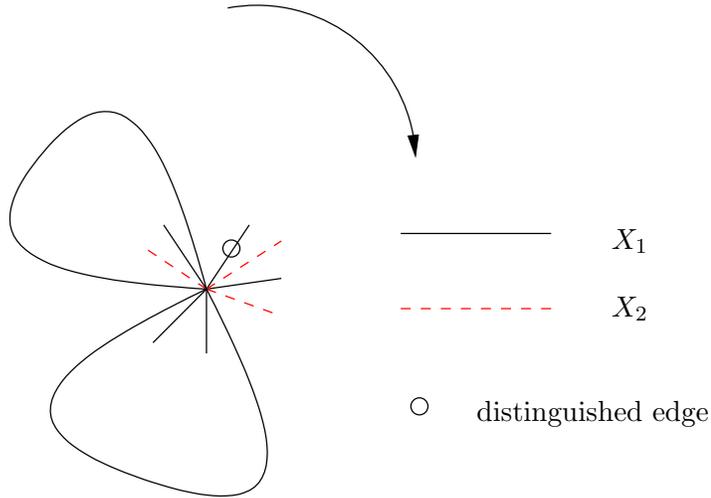}}
\end{center}
\caption{Root of type $X_1X_2X_1X_2\otimes X_1^2\otimes X_2X_1$.}\label{root}
\end{figure}

Thus we define some more complex vertices which will appear 
when we will try to decompose our maps.
Let $P_1$,\dots,$P_\ell$ be a family of monomials.
We associate to this family a bunch of $\ell-1$ circles 
such that outside the circles we put the half-edges of $P_1$ and 
in the $m$-th circle we put the half-edges of a star of type $P_{m+1}$
going out of the central
point and in the same order.
This object will be called the root and we will name each $P_r$
a vertex of the root (look at figure \ref{root}, to see a root of type
$X_1X_2X_1X_2\otimes X_1^2\otimes X_2X_1$).

This corresponds to a star coming from a vertex which have $\ell$ prescribed
loops, the $m$-th having the germs of edges corresponding to a star 
of type $P_r$. Now we construct maps with a root of type
$P_1\otimes\cdots\otimes P_\ell$
and some other vertices of type $q_i$. We say that such a map is minimal if
when we cut the
surface along the $\ell-1$ loops of the root, we do not obtain any component
homeomorphic
to a disc. This means that for any $P_i$ the component of $P_i$ is not planar
i.e. either it
is linked to another $P_j$ or it is linked to some other vertices in a way
that can't be embedded
on a sphere.
\begin{defi}
For $\ell,g$ in $\NN$, a family of monomials $P_1,\cdots,P_\ell$
and integers $\bk=(k_1,\cdots,k_n)$, let
$\MM_g^{\bk}(P_1\otimes\cdots\otimes P_\ell)$ be the number of 
minimal maps of genus $g$ with a root
of type $P_1\otimes\cdots\otimes P_\ell$ and for all $i$, $k_i$ stars of
type $q_i$.
\end{defi}
For example for $V=tX_1X_2X_1X_2$,
the figure \ref{tore} shows that $\MM^1_1(X_1X_2X_1X_2)$ the number of minimal
maps
of genus $1$ with a star of type $X_1X_2X_1X_2$ and a root of type
$X_1X_2X_1X_2$
is $4$.

We extend by linearity $\MM_{\bk}$ and 
$\MM_g^{\bk}$ so that we can compute them on polynomials $P_i$
instead of monomials
and we define
the power series for these enumerations:
$$\I(P)=\sum_{\bk\in\NN^n}
\prod_{i=1}^n\frac{(-t_i)^{k_i}}{k_i!}\MM^{\bk}(P)$$
and
$$\I_g(P_1\otimes\cdots\otimes P_\ell)=\sum_{\bk\in\NN^n}
\prod_{i=1}^n\frac{(-t_i)^{k_i}}{k_i!}
\MM^{\bk}_g(P_1\otimes\cdots\otimes P_\ell).$$
By convention we define for $\l$ in $\CC$, $\I_g(\l)=\l{\mathbbm 1}_{g=0}$ and
$\I_g\equiv 0$ ig $g<0$.

Recall that it was proved in \cite{GMa} that for $\bt$ sufficiently small
$\I(P)=\mu(P)$ for all $P$. This was proved using the fact that these
two quantities satisfy the same induction relation. The induction
relation for the enumeration of maps where given by a decomposition of maps
following the strategy of Tutte. We now try to generalize this fact and we
begin by
looking at the relation given by
decomposing maps.
First, some values can be directly computed
$$\MM_g^{\bk}(1\otimes P_2\otimes\cdots\otimes P_{\ell})
={\mathbbm 1}_{g=\bk=\ell=0}$$
because the component of $1$ is automatically planar.

We now want to count maps that contribute to
$\MM_g^{\bk}(X_iP_1\otimes\cdots\otimes P_{\ell})$ with $P_i$ monomials.
We look at the first half-edge 
of the root $X_iP_1\otimes\cdots\otimes P_{\ell}$ and see where it is glued.
Remember that it must not be planar.

\begin{figure}[ht!]
\psfrag{plus}{$+$}
\psfrag{egal}{''$=$''}
\psfrag{c1}{Where can we glue the first}
\psfrag{c12}{half-edge ?}
\psfrag{c2}{$(1)$: To an other vertex.}
\psfrag{c3}{$(2)$: To a half-edge of the same}
\psfrag{c32}{vertex.}
\psfrag{c4}{$(3)$: To another vertex of}
\psfrag{c42}{the root.}
\psfrag{x1}{$X_1$}
\psfrag{x2}{$X_2$}
\psfrag{de}{distinguished edge}
\begin{center} {\includegraphics{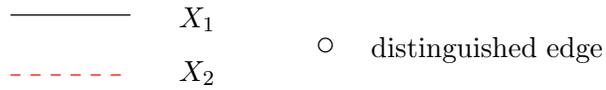}}
\end{center}
\caption{The decomposition process for maps.}\label{dec2}
\end{figure}

Then three cases may occur (see figure \ref{dec2}):
\begin{enumerate}
\item
Either (upper right picture in fig \ref{dec2})
 the half-edge is glued to a vertex of
type $q_j=RX_iS$ for a given $j$.
First we have to choose between the $k_j$ vertices of this
type, then we contract 
the edge coming from this gluing to form a vertex of type 
$SRP_1$. This creates
$$\sum_{1\leq j\leq n, k_j\neq 0}k_j\MM_g^{k_1,\cdots,k_j-1,\cdots,k_n}
(\cycl_iq_jP_1
\otimes P_2\otimes\cdots\otimes P_{\ell})$$
possibilities.
\item
The second case (bottom left picture in fig \ref{dec2}) occurs if the
half-edge is glued
to another half-edge of $P_1=RX_iS$.
It cuts $P_1$ in two: $R$ and $S$. It occurs for all decomposition
of $P_1$ into $P_1=RX_iS$. To write the expression that will arise in
a more convenient way, we will use the 
non-commutative derivative $\part$ which satisfy for $P$ a monomial
$$\part_iP=\sum_{P=RX_iS}R\otimes S.$$
We are now left with two separate circles, one for $R$ and one for $S$.
Either both are non-planar which leads to
$$\sum_{P_1=RX_iS}\MM_g^{\bk}(R\otimes S\otimes P_2
\otimes\cdots\otimes P_{\ell})=\MM_g^{\bk}(\part_iP_1\otimes P_2
\otimes\cdots\otimes P_{\ell})$$
possibilities 
 or one of the component is planar then the two components
can not be linked thus we have to share the vertices of type $q_i$
between them. They are $\left(\sur{\bk}{\bk'}\right)
=\prod_j\left(\sur{k_j}{k'_j}\right)$
ways of choosing for all $j$, $k'_j$ vertices of type $q_j$
for the component of $R$.
 
If the component of $R$ is planar and the one of $S$ is not this
leads to
$$\sum_{{\bk}'+{\bk}''=\bk}\left(\sur{\bk}{\bk'}\right)
\MM_g^{\bk''}(
(\MM^{\bk'}\otimes I)(\part_iP_1)
\otimes\cdots\otimes P_{\ell})$$
possibilities
or $S$ is planar and $R$ is not,
$$\sum_{{\bk}'+{\bk}''=\bk}\left(\sur{\bk}{\bk'}\right)
\MM_g^{\bk'}((I\otimes
\MM^{\bk''})(\part_iP_1)
\otimes\cdots\otimes P_{\ell})$$
possibilities.
\item
The last case occurs if the half-edge is glued with another
vertex $P_r=RX_iS$
of the root. This create a vertex of type $\cycl_iP_rP_1$. Note
that the edge can not cross the circles of the root so it must go
through a handle of the surface thus
it changes the genus by one.
But this vertex is now free from the condition of non planarity so
it can either be planar and thus be separate from the other
vertices of the root:
$$
\sum_{2\leq m\leq l,{\bk}'+{\bk}''=\bk}\left(\sur{\bk}{\bk'}\right)
\MM^{\bk'}(\cycl_iP_rP_1)
\MM_{g-1}^{\bk''}(P_2\otimes\cdots\check{P_r}\cdots\otimes P_{\ell})$$
possibilities
or it may still be non-planar:
$$\sum_{2\leq m\leq \ell}\MM_{g-1}^{\bk}
(\cycl_iP_rP_1\otimes P_2\otimes\cdots\check{P_r}
\cdots\otimes P_{\ell})$$
possibilities.
\end{enumerate}
We can sum these identities and then sum on the $k_i$'s to obtain
the following equality, for all $g$, for all $P_1$\dots,$P_{\ell}$,
\begin{align*}
\I_g(X_iP_1\otimes \cdots\otimes P_{\ell})&=
\sum_j (-t_j)\I_g(\cycl_iq_jP_1\otimes \cdots\otimes P_{\ell})\\
&+\I_g((I\otimes I+\I\otimes I+I\otimes \I)\part_iP_1\otimes \cdots
\otimes P_{\ell})\\
&+\sum_{m\geq 2}\I_{g-1}((I+\I)\cycl_iP_rP_1\otimes\cdots \check{P_r}\cdots
\otimes P_{\ell})
\end{align*}
We can reformulate this by applying it to $P_1=\cycl_i{\bar P}$ and then
summing
on $i$:
\begin{align}\label{relcombinatoire}
\I_g(\Xi P\otimes \cdots\otimes P_{\ell})&=
\sum_{m\geq 2,i}\mu(\cycl_iP_r\cycl_i
{\bar P})\I_{g-1}(P_2\otimes\cdots \check{P_r}\cdots\otimes P_{\ell})
\\
&+\sum_{m\geq 2,i}\I_{g-1}(\cycl_iP_r\cycl_i
{\bar P}\otimes\cdots \check{P_r}\cdots\otimes P_{\ell})\nonumber\\
&+\sum_i\I_g(\part_i\cycl_i{\bar P}\otimes \cdots\otimes P_{\ell}).
\nonumber
\end{align}
where we used the identity $\I=\mu$.
Note that maps that appear in the enumeration must satisfy the condition of 
non-planarity, this imposes a high genus. We have to break the ``planarity''
of $\ell$ components, this can't be done without at least
$\lbrack \frac{\ell+1}{2}\rbrack$
handles on the surface (each handle allow one edge to cross
from one vertex of the root to another one, breaking the planarity
of two components at most).
Thus if $g<Ent(\frac{\ell+1}{2})$,
$\I_g(P_1\otimes \cdots\otimes P_{\ell})=0$. This allow us to 
write the previous equation in a special case which will appear to be useful,
if $\ell$ is even,
\begin{equation}\I_{\frac{\ell}{2}}(\Xi P\otimes \cdots\otimes P_{\ell})=
\sum_{m\geq 2,i}\mu(\cycl_iP_r\cycl_i{\bar P})\I_{\frac{\ell-2}{2}}(P_2\otimes
\cdots \check{P_r}\cdots\otimes P_{\ell}).
\label{relcombeven}
\end{equation}
Thus for $\ell$ even, $\I_{\frac{\ell}{2}}$ satisfy the limit equation
(\ref{limiteven})
of the matrix model. One can easily check that for $\ell$ odd
$\I_{\frac{\hell}{2}}=\I_{\frac{\ell+1}{2}}$ 
satisfy also the limit equation (\ref{limitodd}):
\begin{align*}
\I_{\frac{\ell+1}{2}}(\Xi P\otimes \cdots\otimes P_{\ell})&=
\sum_{m\geq 2,i}\I(\cycl_iP_r\cycl_i\oo{P})\I_{\frac{(\ell-2)+1}{2}}
(P_2\otimes
\cdots \check{P_r}\cdots\otimes P_{\ell})\\
&+\sum_{m\geq 2,i}\I_{\frac{\ell-1}{2}}(\cycl_iP_r\cycl_i\oo{P}\otimes
P_2\otimes
\cdots \check{P_r}\cdots\otimes P_{\ell})\\
&+\sum_i \I_{\frac{\ell+1}{2}}(\partial_i \cycl_i{\bar P}\otimes \cdots\otimes P_{\ell}).
\end{align*}

We can deduce from these identities a control on these enumerations:
\begin{lem}\label{supportcomb}
For all $g\geq 0$, 
there exists $\e>0$ such that for $\bt\in B(0,\e)$,
$\I^g$ is absolutely
convergent and has a bounded support i.e. there exists $M>0$
such that for all polynomials $P_1,\cdots,P_{\ell}$,
$$|\I_g(P_1\otimes \cdots\otimes P_{\ell})|\ppq
\|P_1\|_M\cdots\|P_\ell\|_M$$.
\end{lem}
\begin{dem}
It is sufficient to show that
for all $g\geq 0$, there exists $A_g,B_g>0$ such that for all $h$, for all
monomials $P_1$,\dots,$P_{\ell}$, and all integers
$k_i$:
$$\frac{\MM_g^{\bk}(P_1\otimes\cdots\otimes P_{\ell})}{\prod_i k_i!}
\leq A_g^{\Sigma deg P_i}B_g^{\Sigma k_i}$$
This is easy by induction using the decomposition of maps.
\end{dem}

Finally we need to know the effect of derivation on
these generating function. In fact, derivation adds some vertices
to the enumeration.
\begin{lem}\label{dercomb}
For all $\bj=(j_1,\cdots,j_n)$,
$$\Da_\bj
\I(P)=(-1)^{\bj}\sum_{\bk\in\NN^n}
\prod_{i=1}^n\frac{(-t_i)^{k_i}}{k_i!}\MM^{\bk+\bj}(P)$$
and
$$\Da_\bj
\I_g(P_1\otimes\cdots\otimes P_\ell)=(-1)^{\bj}\sum_{\bk\in\NN^n}
\prod_{i=1}^n\frac{(-t_i)^{k_i}}{k_i!}
\MM^{\bk+\bj}_g(P_1\otimes\cdots\otimes P_\ell).$$
Besides, these series are absolutely convergent and has a bounded support.
\end{lem}
\begin{dem}
The proof is straightforward, $\I$,$\I_g$ are analytic in a neightbourhood
of the origin thus their
derivatives are analytic and their series are given by differentiating
term by term $\I$ and $\I_g$.
\end{dem}
Thus derivatives fix some vertices in the enumeration a fact often used
in combinatorics to find relation between generating functions of graphs.

\section{High order observable}
We have already seen that
$\nu(P_1\otimes\cdots\otimes P_\ell)=\I_{\frac{\hell}{2}}
(P_1\otimes\cdots\otimes P_\ell)$ satisfy the
limit equation of $\nu^N$. This is our candidate for the limit of the $\nu^N$'s.
In fact this suggests a statement closely related to \ref{theofree}.
\begin{prop}\label{mainprop}
For all $\ell$, for all $g\in \NN$, there exists $\eta>0$
such that for all $\bt$ in
$B_{\eta,c}$, for all polynomials $P_1$,\dots,$P_\ell$,
\begin{align*}
\nu^N(P_1\otimes\cdots\otimes P_\ell)=
\I_{\frac{\hell}{2}}(P_1\otimes\cdots\otimes P_\ell)
+\frac{1}{N^2}\I_{\frac{\hell}{2}+1}(P_1\otimes\cdots\otimes P_\ell)+\\
\cdots
+\frac{1}{N^{2g}}\I_{\frac{\hell}{2}+g}(P_1\otimes\cdots\otimes P_\ell)
+o(\frac{1}{N^{2g}}).
\end{align*}
\end{prop}

To prove this we have to define all the correction to the convergence.
We define $\nu^N_1=\nu^N$ and by induction on $h$, for all $N$,
all polynomials
$P_1,\dots,P_\ell$,
$$\nu_{h+1}^{N}(P_1\otimes\cdots\otimes P_\ell)
=N^2(\nu_{h}^{N}
-\I_{\frac{\hell}{2}+h-1})(P_1\otimes\cdots\otimes P_\ell).$$
Those quantities satisfy also some induction relation similar to those
of property \ref{sd1}
\begin{prop}\label{sdg}
For all $h\pgq 2$, $\ell$ in $\NN$, for all polynomial
$P_1$,\dots,$P_\ell$, for all $N$, the ``finite Schwynger-Dyson's
equation of order $h$''  $({\bf SD}^{N}_{h,\ell})$ is satisfied by $\nu^N$:
if $\ell$ is even
\begin{align*}\label{releven}
\nu^N_h(\Xi P\otimes P_2\otimes \cdots\otimes P_\ell)
&=\sum_{i,r}\mu(\cycl_i\oo{P}\cycl_i P_r)\nu^N_h
( P_2\otimes \cdots \check{P_r} \cdots\otimes P_\ell)\\
&\hspace{-3cm}+\sum_{i,r}\nu^N_{h-1}
(\cycl_i\oo{P}\cycl_i  P_r\otimes P_2\otimes \cdots \check{P_r} \cdots\otimes
P_\ell)\hspace{2.5cm}({\bf SD}^{N}_{h,\ell})\\
&\hspace{-3cm}+\sum_i\nu^N_{h-1}
(\part_i\cycl_i\oo{P}\otimes P_2\otimes\cdots\otimes P_\ell)
\end{align*}
and if $\ell$ is odd
\begin{align*}
\nu^N_h(\Xi P\otimes P_2\otimes \cdots\otimes P_\ell)
&=\sum_{i,r}\mu(\cycl_i\oo{P}\cycl_i P_r)\nu^N_h
( P_2\otimes \cdots \check{P_r} \cdots\otimes P_\ell)\\
&\hspace{-3cm}+\sum_{i,r}\nu^N_h
(\cycl_i\oo{P}\cycl_i  P_r\otimes P_2\otimes \cdots \check{P_r} \cdots\otimes
P_\ell)\hspace{2.5cm} ({\bf SD}^{N}_{h,\ell})\\
&\hspace{-3cm}+\sum_i\nu^N_h
(\part_i\cycl_i\oo{P}\otimes P_2\otimes\cdots\otimes P_\ell)
\end{align*}
\end{prop}
\begin{dem}
Remember that we have shown in Property \ref{sd1}, for $\ell$ even
\begin{align*}
\nu^N_1(\Xi P\otimes P_2\otimes \cdots\otimes P_\ell)
&=\sum_{i,r}\mu_{\bt}(\cycl_i\oo{P}\cycl_i P_r)\nu^N_1
( P_2\otimes \cdots \check{P_r} \cdots\otimes P_\ell)\\
&+\frac{1}{N^2}\sum_{i,r}\nu^N_1
(\cycl_i\oo{P}\cycl_i  P_r\otimes P_2\otimes \cdots \check{P_r} \cdots\otimes
P_\ell)\\
&+\frac{1}{N^2}\sum_i\nu^N_1
(\part_i\cycl_i\oo{P}\otimes P_2\otimes\cdots\otimes P_\ell)
\end{align*}
and according to (\ref{relcombeven})
$$\I_{\frac{\ell}{2}}(\Xi P\otimes \cdots\otimes P_{\ell})=
\sum_{m\geq 2,i}\I(\cycl_iP_r\cycl_i{\bar P})\I_{\frac{\ell-2}{2}}(P_2\otimes
\cdots \check{P_r}\cdots\otimes P_{\ell}).$$
Thus if we substract these two equalities and multiply the result by
$N^2$ we obtain (${\bf SD}^{N}_{2,\ell}$) (Observe that with our convention
$\nu^N(\l)=\I_0(\l)$).
\begin{align*}
\nu^N_2(\Xi P\otimes P_2\otimes \cdots\otimes P_\ell)
&=\sum_{i,r}\mu(\cycl_i\oo{P}\cycl_i P_r)\nu^N_2
( P_2\otimes \cdots \check{P_r} \cdots\otimes P_\ell)\\
&+\sum_{i,r}\nu^N_{1}
(\cycl_i\oo{P}\cycl_i  P_r\otimes P_2\otimes \cdots \check{P_r} \cdots\otimes
P_\ell)\\
&+\sum_i\nu^N_{1}
(\part_i\cycl_i\oo{P}\otimes P_2\otimes\cdots\otimes P_\ell).
\end{align*}

Now suppose that for $\ell$ even, $h\pgq 2$, for all polynomial
$P_1$,\dots,$P_\ell$, for all $N$,
(${\bf SD}^{N}_{h,\ell}$) is satisfied.
Then according to (\ref{relcombinatoire}):
\begin{align*}
\I_{\frac{\ell}{2}+h-1}(\Xi P\otimes \cdots\otimes P_{\ell})&=
\sum_i\I_{\frac{(\ell+1)+1}{2}+h-2}(\part_i \cycl_i{\bar P}\otimes
\cdots\otimes P_{\ell})\\
&\hspace{-3cm}+\sum_{r\geq 2,i}\I_{\frac{(\ell-1)+1}{2}+h-2}
(\cycl_iP_r\cycl_i{\bar P}
\otimes\cdots \check{P_r}\cdots\otimes P_{\ell})\\
&\hspace{-3cm}+\sum_{r\geq 2,i}\mu(\cycl_iP_r\cycl_i{\bar P})
\I_{\frac{(\ell-2)}{2}+h-1}(\cycl_iP_r\cycl_i
{\bar P}\otimes\cdots \check{P_r}\cdots\otimes P_{\ell}).
\end{align*}
and this can be translated into
\begin{align*}
\I_{\frac{\hat{\ell}}{2}+h-1}(\Xi P\otimes \cdots\otimes P_{\ell})&=
\sum_i\I_{\frac{\widehat{\ell+1}}{2}+h-2}(\part_i \cycl_i{\bar P}\otimes
\cdots\otimes P_{\ell})\\
&\hspace{-3cm}+\sum_{r\geq 2,i}\I_{\frac{\widehat{\ell-1}}{2}+h-2}
(\cycl_iP_r\cycl_i{\bar P}
\otimes\cdots \check{P_r}\cdots\otimes P_{\ell})\\
&\hspace{-3cm}+\sum_{r\geq 2,i}\mu(\cycl_iP_r\cycl_i{\bar P})
\I_{\frac{\widehat{\ell-2}}{2}+h-1}(\cycl_iP_r\cycl_i
{\bar P}\otimes\cdots \check{P_r}\cdots\otimes P_{\ell}).
\end{align*}
Substracting this equality from (${\bf SD}^{N}_{h,\ell}$) we get
(${\bf SD}^{N}_{h+1,\ell}$).
This proves by induction (${\bf SD}^{N}_{h,\ell}$) for all $h$, and for
all $\ell$ even.

We proceed in the same way for $\ell$ odd. Observe that the equation for
$\ell$ odd and $h=1$ is satisfied
according to Property \ref{sd1}. Then observe that for $\ell$ odd,
$({\bf SD}^{N}_{h+1,\ell})$ can be obtained by substracting
(\ref{relcombinatoire}) with $g=\frac{\ell+1}{2}+h-1$
\begin{align*}
\I_{\frac{\hat{\ell}}{2}+h-1}(\Xi P\otimes \cdots\otimes P_{\ell})&=
\sum_i\I_{\frac{\widehat{\ell+1}}{2}+h-2}(\part_i \cycl_i{\bar P}\otimes
\cdots\otimes P_{\ell})\\
&\hspace{-3cm}+\sum_{r\geq 2,i}\I_{\frac{\widehat{\ell-1}}{2}+h-2}
(\cycl_iP_r\cycl_i{\bar P}
\otimes\cdots \check{P_r}\cdots\otimes P_{\ell})\\
&\hspace{-3cm}+\sum_{r\geq 2,i}\mu(\cycl_iP_r\cycl_i{\bar P})
\I_{\frac{\widehat{\ell-2}}{2}+h-1}(\cycl_iP_r\cycl_i
{\bar P}\otimes\cdots \check{P_r}\cdots\otimes P_{\ell}).
\end{align*}
 from
(${\bf SD}^{N}_{h,\ell}$).
\end{dem}

\section{Asymptotic of the matrix model}
The issue with the previous relations is that
they only give us the  moments of products of polynomials
such that the first polynomial is in the image of $\Xi$. Thus
we need to invert $\Xi$.
We define the operator norm with respect to $\|.\|_M$:
$$|||A|||_M=\sup_{\|P\|_M\ppq 1}\|AP\|_M.$$
In \cite{GMa}, we give some estimates on the operator norm of $\Xi$.
\begin{lem}\label{decrease}
\begin{enumerate}
\item
The operator $\Xi_0$ is invertible on $\cxm$.
\item
There exists $M_0>0$ such that for all $M>M_0$, the operators 
$\Xi_2$, $\Xi_0$
and $\Xi_0^{-1}$ are continuous
and their norm are uniformly 
bounded for $\bt$ in $B_{\eta}$.
\item
For all polynomials $P$, $\deg \Xi_0^{-1}P\ppq \deg P$ and
$\deg \Xi_1 P\ppq \deg P+\deg V-2$
\item
For all $\e,M>0$, there exists $\eta_{\e}>0$ such for $|\bt|<\eta_{\e}$,
$\Xi_1$ is continuous
on $\cxm$  and
$|||\Xi_1|||_{M}\ppq \e$.
\end{enumerate}
\end{lem}

The last step to proves Theorem \ref{theofree} is to control the
$\nu_N^h$.

\begin{lem}\label{sharpestimate}
For all $\ell,h\in\NN^*$, $\a>0$, there exists $C,\eta,M_0>0$,
such that for all $\bt\in B_{\eta,c}$, $M>M_0$ and
all polynomials $P_1$,\dots,$P_\ell$ of degree less than $\a N^{\frac{1}{2}}$,
$$|\nu^N_h(P_1\otimes\cdots\otimes P_\ell)|\leq C\|P_1\|_M\cdots\|P_\ell\|_M$$
\end{lem}
\begin{dem}
The case $h=1$, $\ell$ even is a direct consequence of Lemma \ref{estimate}.
We treat the other cases by induction using Property \ref{sdg}.
As the equations are different according to the parity of $\ell$, 
we have to be careful: we prove the result by an induction on $h$ and for
a fixed
$h$ we deal first with the case $\ell$ even and then with the case $\ell$
odd (Note that
both time we will do an induction on $\ell$).
Now we choose $\ell,h\in\NN^*$, $\a>0$ and polynomials
$P_1$,\dots,$P_\ell$ of degree less than $\a N^{\frac{1}{2}}$.

Then, the idea to nearly invert $\Xi$ on a polynomial $P$
is to approximate $P_1$ by
$\Xi Q_n=(\Xi_0+\Xi_1)Q_n$ with
$$Q_n=\sum_{k=0}^{n-1} (-{ \Xi}_0^{-1}{ \Xi}_1)^k{ \Xi}_0^{-1} P_1.$$
The remainder is
$$R_n = P_1- \Xi Q_n= (-{ \Xi}_1{ \Xi}_0^{-1})^{n} P_1.$$
As $\Xi_1$ is the multiplication by a derivative of $V$ it should have 
a small norm and the remainder should be easily controlled.

We can make the decomposition:
\begin{equation}\label{dec}\nu_h^N(P_1\otimes \cdots\otimes P_\ell)
=\nu_h^N(\Xi Q_n\otimes \cdots\otimes P_\ell)
+\nu_h^N(R_n\otimes \cdots\otimes P_\ell)
\end{equation}
Now we let $n$ goes to infinity with $N$, for example $n=[\sqrt{N}]$.
It is important that $n$ goes to infinity no too slowly but we must have
$n=O(\sqrt{N})$ in order to use all the induction hypothesis.
An important fact is that the degrees of $R_n$ and $Q_n$ are
 $O(\sqrt{N})$ since
$\Xi_0^{-1}\Xi_1$ change the degree by at most $D-2$.

We first control the term with $R_n$, by definition of the $\nu_N^h$,
$$\nu_h^N(R_n\otimes \cdots\otimes P_\ell)=(N^{2(h-1)}\nu_1^N
+N^{2(h-1)}\I_{\frac{\hat{\ell}}{2}}+\dots+
\I_{\frac{\hat{\ell}}{2}
+h-1})
(R_n\otimes \cdots\otimes P_\ell).
$$
Each of the $\I_g$ are compactly supported according to
Lemma \ref{supportcomb} so that if $\eta$  is sufficiently small,
for $\bt$ in $B_{\eta,c}$, $\I_{\frac{\hat{\ell}}{2}}$,
\dots,$\I_{\frac{\hat{\ell}}{2}+h-1}$ are convergent and we can take
$M$ bigger than the radius of their support.
Besides, Property \ref{estimate} shows that for polynomials $P$ of degree
of order $N^{\frac{1}{2}}$,
$$|\nu_1^N(P_1\otimes \cdots\otimes P_\ell)|\ppq CN\|P_1\|_M \cdots
\|P_\ell\|_M.$$
Thus according to Lemma \ref{decrease}, for
$\eta$ small, $\|R_n\|_M\ppq|||{ \Xi}_1{ \Xi}_0^{-1}|||^n\|P_1\|_M$
decrease exponentially fast in $n$ and this uniformly for $\bt\in B_{\eta,c}$.
Then since $n\sim\sqrt{N}$
$$|\I_g(R_n\otimes \cdots\otimes P_\ell)|
\ppq \|R_n\|_M\|P_2\|_M\dots\|P_\ell\|_M
\ppq Ce^{-C'\sqrt{N}}\|P_1\|_M\dots\|P_\ell\|_M.$$
Thus,
\begin{equation}\label{remainder}
|\nu_h^N(R_n\otimes \cdots\otimes P_\ell)|
\ppq N^{2h}Ce^{-C'\sqrt{N}}\|P_1\|_M\dots\|P_\ell\|_M
\end{equation}
and $N^{2h}Ce^{-C'\sqrt{N}}$ is bounded.

Finally we have to deal with $\nu_h^N(\Xi Q_n\otimes \cdots\otimes P_\ell)$. 
We can use $({\bf SD}_{h,\ell}^N)$:
\begin{align}\label{rec}
\nu^N_h(\Xi Q_n\otimes P_2\otimes \cdots\otimes P_\ell)
&=\sum_{i,m}\mu_{\bt}(\cycl_i{\bar Q_n}\cycl_i P_r)\nu^N_{h}
( P_2\otimes \cdots \check{P_r} \cdots\otimes P_\ell)\\
&+\sum_{i,m}\nu^N_{h-{{\mathbbm 1}_{\ell\textrm{ even}}}}
(\cycl_i{\bar Q_n}\cycl_i  P_r\otimes P_2\otimes \cdots
\check{P_r} \cdots\otimes P_\ell)\nonumber\\
&
+\sum_i\nu^N_{h-{{\mathbbm 1}_{\ell\textrm{ even}}}}
(\part_i\cycl_i{\bar Q_n}\otimes P_2\otimes\cdots\otimes P_\ell)\nonumber
\end{align}

We now use the induction hypothesis. Indeed if $h$ is even, on the
right hand side either $h$ decreases or $h$ remains constant and $\ell$
decreases and remains even. If $h$ is odd, either $\ell$ becomes even
or $\ell$ decreases.
Now, let $C$ be an uniform bound on the norm of $\Xi_0^{-1}$ (which
exists according to Lemma \ref{decrease}) by definition of $Q_n$,
$$\|Q_n\|_M\ppq\sum_{k=0}^{n-1}
|||{ \Xi}_0^{-1}{ \Xi}_1|||_M^k|||{ \Xi}_0^{-1}|||_M\| P\|_M
\ppq\frac{C}{1-C|||\Xi_1|||_M}\|P\|_M.$$
Thus, using Lemma \ref{decrease}, if $\eta$ is sufficiently small,
for $\bt$ in $B_{\eta,c}$, $\|Q_n\|_M\ppq 2\|P\|_M$.
Note that
$$\deg Q_n\ppq \deg P_1+2\sqrt{N}(D-1)\ppq(\alpha+2(D-1))\sqrt{N}.$$
We can now apply the induction hypothesis with $\alpha'=\alpha+2(D-1)$,
there exists $M,C,\eta$ such that for $\bt$ in $B_{\eta,c}$,
\begin{align}\label{ineq1}
&|\mu_{\bt}(\cycl_i{\bar Q_n}\cycl_i P_r)\nu^N_{h}
( P_2\otimes \cdots \check{P_r} \cdots\otimes P_\ell)|\\
&\ppq C\|\cycl_i{\bar Q_n}\cycl_i P_r\|_M\|P_2\|_M\cdots
\|\check{P_r}\|_M\cdots  \|P_{\ell}\|\nonumber
\end{align}
where we have assumed that $M$ is bigger than the radius of the support
of $\mu_{\bt}$ which is bounded according to (\ref{support}).
Besides, by the induction hypothesis we can obtain with
the same constant,
\begin{align}\label{ineq2}
&|\nu^N_{h-{{\mathbbm 1}_{\ell\textrm{ even}}}}
(\cycl_i{\bar Q_n}\cycl_i  P_r\otimes P_2\otimes \cdots
\check{P_r} \cdots\otimes P_\ell)|\\
&\ppq C\|\cycl_i{\bar Q_n}\cycl_i P_r\|_M\|P_2\|_M\cdots
\|\check{P_r}\|_M\cdots  \|P_{\ell}\|_M\nonumber
\end{align}
and
\begin{align}\label{ineq3}
&|\nu^N_{h-{{\mathbbm 1}_{\ell\textrm{ even}}}}
(\part_i\cycl_i{\bar Q_n}\otimes P_2\otimes\cdots\otimes P_\ell)|\\
&\ppq \|\part_i\cycl_i{\bar Q_n}\|_M\|P_2\|_M\cdots
\|\check{P_r}\|_M\cdots \|P_{\ell}\|_M.\nonumber
\end{align}
Now remember, that if $M'>M$,
$$\cycl_i:(\cxm,\|.\|_{M'})\to(\cxm,\|.\|_M)$$
is continuous, thus
\begin{align*}\|\cycl_i{\bar Q_n}\cycl_i P_r\|_M
&\ppq \|\cycl_i{\bar Q_n}\|_M\|\cycl_i P_r\|_M\\
&\ppq C\|Q_n\|_{M'}\|P_r\|_{M'}\ppq C\|P_1\|_{M'}\|P_r\|_{M'}
\end{align*}
and
$$\|\part_i\cycl_i {\bar Q_n}\|_M\ppq C\|P_1\|_{M'}.$$
If we use inequalities (\ref{ineq1}), (\ref{ineq2}) and (\ref{ineq3})
in the decomposition (\ref{rec}), we get
\begin{equation*}
|\nu^N_h(\Xi Q_n\otimes P_2\otimes \cdots\otimes P_\ell)|
\ppq C\|P_1\|_{M'}\|P_2\|_{M'}\cdots
 \|P_{\ell}\|_{M'}.
\end{equation*}
Finally, we conclude with (\ref{remainder}) and (\ref{dec})
\end{dem}

Now, for all $\ell$, $h$ there exists $\eta>0$ such that for
$\bt\in B_{\eta,c}$, for all polynomials $P_1,\cdots,P_\ell$,
$$\nu^N_g(P_1\otimes\cdots\otimes P_{\ell})
=N^2(\nu^N_g-\I_{\frac{\hat{\ell}}{2}+g-1})
(P_1\otimes\cdots\otimes P_{\ell})$$
is a bounded sequence for all $g\ppq h+1$. Thus for all $g\ppq h$,
$\nu^N_g(P_1\otimes\cdots\otimes P_{\ell})$ goes to
$\I_{\frac{\hat{\ell}}{2}+g-1}$ and
$$\nu^N(P_1\otimes\cdots\otimes P_{\ell})=
\I_{\frac{\hat{\ell}}{2}}
+\frac{1}{N^2}\I_{\frac{\hat{\ell}}{2}+1}+\cdots+
+\frac{1}{N^{2h}}\I_{\frac{\hat{\ell}}{2}+h}+o(\frac{1}{N^{2h}}).$$
Thus Property \ref{mainprop} is proved. The special case $\ell=1$ is
exactly Theorem \ref{maintheo}:
\begin{align*}
E[\mun(P)]&=\I(P)+\frac{1}{N^2}\nu^N_1(P)\\
&=\I(P)+\frac{1}{N^2}\I_1(P)+\cdots+\frac{1}{N^{2g}}\I_g(P)
+o(\frac{1}{N^{2g}}).
\end{align*}

Thus we can prove Theorem \ref{theofree}.
\begin{theo}
For all $g\in \NN$, there exists $\eta>0$
such that for all $\bt$ in
$B_{\eta,c}$,
$$F^N_{V_{\bt}}
=F^0(\bt)+\cdots
+\frac{1}{N^{2g}}F^g(\bt)+o(\frac{1}{N^{2g}})$$
and $F^g$ is
the generating function for
maps of genus $g$ associated with $V$:
$$F^g(\bt)=
\sum_{k_1,\cdots,k_n\in\NN}
\prod_i\frac{(-t_i)^{k_i}}{k_i!}
\C^{\bk}$$
where if $\bk=(k_1,\cdots,k_n)$, $\C^{\bk}$ is the number of maps
on a surface
of
genus $g$ with
$k_i$ vertices of type $q_i$.
\end{theo}
\begin{dem}
Note that the estimate we get in Lemma \ref{derestimate} are uniform
in $\bt$ provided we are in $B_{\eta,c}$. Now observe that if $V_{\bt}$
is $c$-convex then for $\alpha$ in $[0,1]$,
$V_{\alpha \bt}$ is $c$-convex if $c\ppq 1$ and $1$-convex if $c>1$.
Thus if $\bt$ is in $B_{\eta,c}$, for all $0<\alpha<1$,
$\alpha\bt\in B_{\eta,\min(c,1)}$. This allow us to use Property \ref{propder}
with an uniformly bounded remainder.
\begin{align*}
F^N_{V_{\bt}}&=
\int_0^1-E_{\mu^N_{V_{\alpha\bt}}}[\mun(V_{\bt})]d\alpha\\
&=-\int_0^1\mu_{\alpha\bt}(V_{\bt})d\alpha-\frac{1}{N^2}\int_0^1
\nu^N_{\alpha\bt}(\mun(V_{\bt}))d\alpha\\
&=F^0(\bt)+\cdots
+\frac{1}{N^{2g}}F^g(\bt)+o(\frac{1}{N^{2g}})
\end{align*}
with
\begin{align*}
F^0&=-\int_0^1
\sum_{k_1,\cdots,k_n\in\NN}
\prod_i\frac{(-\alpha t_i)^{k_i}}{k_i!}
\MM^{\bk}(V_{\bt}) d\alpha\\
&= \sum_{k_1,\cdots,k_n\in\NN}
\prod_i\frac{(-t_i)^{k_i}}{k_i!}
\C^{\bk}
\end{align*}
since it can be easily checked that
$$\sum_{k_1,\cdots,k_n\in\NN}
\prod_i\frac{(-u t_i)^{k_i}}{k_i!}
\C^{\bk}$$ has the same derivative in $u$ than
$-\int_0^u
\sum_{k_1,\cdots,k_n\in\NN}
\prod_i\frac{(-\alpha t_i)^{k_i}}{k_i!}
\MM^{\bk}(V_{\bt}) d\alpha.$
With the same technique,
\begin{align*}
F^g&=-\int_0^1
\sum_{k_1,\cdots,k_n\in\NN}
\prod_i\frac{(-\alpha t_i)^{k_i}}{k_i!}
\MM^{\bk}_g(V_{\bt}) d\alpha\\
&= \sum_{k_1,\cdots,k_n\in\NN}
\prod_i\frac{(-t_i)^{k_i}}{k_i!}
\C^{\bk}_g
\end{align*}
This proves the Theorem:
$$
F^N_{V_{\bt}}=F^0+\cdots+\frac{1}{N^2g}F^g+o(\frac{1}{N^{2g}}).$$
\end{dem}

\section{Higher derivatives.}
In this section we will show that
one can derivate these expansions term by terms. Indeed, the family of
the $\nu^N_h$'s is sufficiently rich to express any of its 
own derivatives. Thus,
we will be able to find a recursive decomposition of this derivatives.

\begin{prop}\label{sdder}
For all $1\ppq j\ppq n$, for all polynomials $P_1,\cdots,P_\ell$,
if $\ell$ is even,
\begin{align*}
\frac{\partial}{\partial t_j}\nu^N_1(P_1\otimes\cdots\otimes P_{\ell})
&=-\nu^N_1(P_1\otimes\cdots\otimes P_{\ell}\otimes q_j)\\
&+\sum_{r=1}^\ell \frac{\partial}{\partial t_j}\mu(P_r)
\nu^N_1(P_1\otimes\cdots\check{P_r}\cdots\otimes P_{\ell})\\
&+\nu^N_1(P_1\otimes\cdots\otimes P_{\ell})\nu^N_1(q_j)
\end{align*}
and if $\ell$ is odd,
\begin{align*}
\frac{\partial}{\partial t_j}\nu^N_1(P_1\otimes\cdots\otimes P_{\ell})
&=-\nu^N_2(P_1\otimes\cdots\otimes P_{\ell}\otimes q_j)\\
&+\sum_{r=1}^\ell \frac{\partial}{\partial t_j}\mu(P_r)
\nu^N_2(P_1\otimes\cdots\check{P_r}\cdots\otimes P_{\ell})\\
&+\nu^N_1(P_1\otimes\cdots\otimes P_{\ell})\nu^N_1(q_j).
\end{align*}
\end{prop}
\begin{dem}
We simply need to derivate 
$$\frac{1}{Z^N_V}\int\prod_r(\frac{1}{N}\tr-\mu)(P_r)e^{-N\tr V}d\mu^N.$$
In that expression we can either derivate $\frac{1}{Z^N_V}$, the potential 
$e^{-N\tr V}$ or one of the $\mu(P_r)$, this leads to
\begin{align*}
\frac{\partial}{\partial t_j}
E[\prod_r(\frac{1}{N}\tr-\mu)(P_r)]&=
N^2E[\prod_r(\frac{1}{N}\tr-\mu)(P_r)(\frac{1}{N}\tr-\mu)(q_j)]\\
&-E[\prod_r(\frac{1}{N}\tr-\mu)(P_r)]N^2E[(\frac{1}{N}\tr-\mu)(q_j)]\\
&+\sum_r \frac{\partial}{\partial t_j}\mu(P_r)
E[\prod_{r'\neq r}(\frac{1}{N}\tr-\mu)(P_{r'})]
\end{align*}
Where one can notice that we have added in the two first terms of the
right hand side the quantity $\mu(q_j)$ but these two modifications
cancel each other.

Now multiply by the normalisation $N^\ell$ to get the equation
in the case $\ell$ even. In the case $\ell$ odd, if we multiply by 
$N^{\ell+1}$ we get
\begin{align*}
\frac{\partial}{\partial t_j}\nu^N_1(P_1\otimes\cdots\otimes P_{\ell})
&=-N^2\nu^N_1(P_1\otimes\cdots\otimes P_{\ell}\otimes q_j)\\
&+N^2\sum_{r=1}^\ell \frac{\partial}{\partial t_j}\mu(P_r)
\nu^N_1(P_1\otimes\cdots\check{P_r}\cdots\otimes P_{\ell})\\
&+\nu^N_1(P_1\otimes\cdots\otimes P_{\ell})\nu^N(q_j)\\
&=-\nu^N_2(P_1\otimes\cdots\otimes P_{\ell}\otimes q_j)\\
&+\sum_{r=1}^\ell \frac{\partial}{\partial t_j}\mu(P_r)
\nu^N_2(P_1\otimes\cdots\check{P_r}\cdots\otimes P_{\ell})\\
&+\nu^N_1(P_1\otimes\cdots\otimes P_{\ell})\nu^N(q_j)\\
&+N^2r_N
\end{align*}
with by definition of $\nu^N_2$,
$$r_N=-\I_{\frac{l+1}{2}}(P_1\otimes\cdots\otimes P_{\ell}\otimes q_j)
+\sum_r \frac{\partial}{\partial t_j}\mu(P_r)
\I_{\frac{l-1}{2}}(P_1\otimes\cdots\check{P_r}\cdots\otimes P_{\ell}).$$
Let's have a closer look at this expression,
$\I_{\frac{l+1}{2}}(P_1\otimes\cdots\otimes P_{\ell}\otimes q_j)$ counts maps
with $\frac{l+1}{2}$ handles and such none of the $l+1$ components 
$P_1$,\dots, $P_\ell$, $q_j$ is planar. Since each handle can break
the planarity of at most two components,
the only way to obtain such a configuration is to form
 $\frac{l+1}{2}$ couples among these $l$ components and in each of these
couples to put a handle between the two components.
For example, one can decompose these maps
according to the other vertex in the couple of $q_j$:
$$\I_{\frac{l+1}{2}}(P_1\otimes\cdots\otimes P_{\ell}\otimes q_j)
=\sum_r \I_1(P_r\otimes q_j)
\I_{\frac{l-1}{2}}(P_1\otimes\cdots\check{P_r}\cdots\otimes P_{\ell}).$$
Now observe that $\I_1(P_r\otimes q_j)$ counts maps of genus $1$ such
that the component of $P_r$ is linked to the component of $q_j$, thus
it is equivalent to the counting of planar maps with two precribed vertices,
one of type $P_r$ and one of type $q_j$. According to Lemma \ref{dercomb},
that exactly what count $\frac{\partial}{\partial t_j}\mu(P_r)$.
Thus $r_N=0$ and the Property is proved.
\end{dem}

With this decomposition we can now show that the derivatives of the
$\nu^N_h$'s are of order $1$.

\begin{lem}\label{derestimate}
For all
$\bj=(j_1,\cdots,j_n)$
, 
for all $\ell,h\in\NN^*$, $\a>0$, there exists constants
$C,\eta,M_0>0$,
such that for all $\bt\in B_{\eta,c}$, $M>M_0$ and
all polynomials $P_1$,\dots,$P_\ell$ of degree less than $\a N^{\frac{1}{2}}$,
$$|\Da_\bj\nu^N_h(P_1\otimes\cdots\otimes P_\ell)|\leq C\|P_1\|_M\cdots
\|P_\ell\|_M$$
\end{lem}
\begin{dem}
The proof is essentially the same than the proof of Lemma \ref{sharpestimate}.
We just use in
addition
an induction on $\sum_i j_i$. If $\sum_i j_i=0$ then we are exactly in the
case of the Lemma \ref{sharpestimate}.
Otherwise, Assume that
we have already proved this lemma up to the $\sum_i j_i-1$-th derivative.
To prove
it for $\sum_i j_i$, we proceed by induction.
First we prove it for $\nu^N_1$, and it is straightforward since we can
express the derivatives of $\nu^N_1$ as a function of derivatives of lesser
degree according to Property \ref{sdder}.
Then for $h>1$ we
need to find an
induction
relation on the
$\Da_\bj\nu_N^h$.
We start form Property \ref{sdg}:
\begin{align*}
\nu^N_h(\Xi P_1\otimes P_2\otimes \cdots\otimes P_\ell)
&=\sum_{i,m}\mu_{\bt}(\cycl_i{\bar P_1}\cycl_i P_r)\nu^N_{h}
( P_2\otimes \cdots \check{P_r} \cdots\otimes P_\ell)\\
&+\sum_{i,m}\nu^N_{h-{{\mathbbm 1}_{\ell\textrm{ even}}}}
(\cycl_i{\bar P_1}\cycl_i  P_r\otimes P_2\otimes \cdots
\check{P_r} \cdots\otimes P_\ell)\\
&
+\sum_i\nu^N_{h-{{\mathbbm 1}_{\ell\textrm{ even}}}}
(\part_i\cycl_i{\bar P_1}\otimes P_2\otimes\cdots\otimes P_\ell)
\end{align*}
and we take the derivative:
\begin{align*}
&\Da_{\bj}\nu^N_h(\Xi P_1\otimes P_2\otimes \cdots\otimes P_\ell)\\
&=\sum_{\sur{\bj'+\bj''=\bj}{\bj''\neq 0}}(\sur{\bj'}{\bj})
\Da_{\bj'}\nu^N_h((\Da_{\bj''}
\Xi) P_1\otimes P_2\otimes \cdots\otimes P_\ell)\\
&+\sum_{\sur{\bj'+\bj''=\bj}{\bj''\neq 0}}\sum_{i,m}(\sur{\bj'}{\bj})
\Da_{\bj'}\mu_{\bt}(\cycl_i{\bar P_1}\cycl_i P_r)\Da_{\bj''}\nu^N_{h}
( P_2\otimes \cdots \check{P_r} \cdots\otimes P_\ell)\\
&+\sum_{i,m}\Da_{\bj}\nu^N_{h-{{\mathbbm 1}_{\ell\textrm{ even}}}}
(\cycl_i{\bar P_1}\cycl_i  P_r\otimes P_2\otimes \cdots
\check{P_r} \cdots\otimes P_\ell)\\
&
+\sum_i\Da_{\bj}\nu^N_{h-{{\mathbbm 1}_{\ell\textrm{ even}}}}
(\part_i\cycl_i{\bar P_1}\otimes P_2\otimes\cdots\otimes P_\ell)
\end{align*}
with $$\Da_{\bj''}\Xi P=(I\otimes(\Da_{\bj''}\mu)+(\Da_{\bj''}\mu)\otimes I)\part_i\cycl_i {\bar P}-(\Da_{\bj''}\cycl_iV)\cycl_i{\bar P}.$$
Thus $\Da_{\bj}\nu^N_h$ has a decomposition similar to the one of
Property \ref{sdg} and we can use it in  similar way we use it in the
proof of Lemma \ref{sharpestimate}. The rest of the proof is identical
to the one of Lemma \ref{sharpestimate} and we only give the main steps.

We define the same $Q_n$
and $R_n$ as in the previous proof to approximate $\Xi^{-1}P$.
Then, $\Da_{\bj}\nu^N_h(\Xi Q_n\otimes P_2\otimes \cdots\otimes P_\ell)$
is controlled using the previous decomposition, using the fact that
$\Da_{\bj}\mu$ has a bounded support and that the norm of
$(\Da_{\bj''}\cycl_iV)$ is also bounded.
Now write 
\begin{align*}
\Da_{\bj}\nu^N_h(R_n\otimes P_2\otimes \cdots\otimes P_\ell)
&=(\Da_{\bj} N^{2(h-1)}\nu^N_1+N^{2(h-1)}\Da_{\bj}\I_{\hat{\ell}}+\cdots\\
&
+\Da_{\bj}\I_{\hat{\ell}+h-1})(R_n\otimes P_2\otimes \cdots\otimes P_\ell)
\end{align*}
Then for $\bt$ sufficiently small these series are convergent and have
uniformly bounded support according to Lemma \ref{dercomb} and we can use
Property \ref{sdder} to control one more time
$\Da_{\bj}\nu^N_1$.
\end{dem}

From there we deduce
\begin{prop}\label{propder}
For all $\ell$, for all $g\in \NN$, for all
$\bj=(j_1,\cdots,j_n)$ there exists $\eta>0$
such that for all $\bt$ in
$B_{\eta,c}$, for all polynomials $P_1$,\dots,$P_\ell$,
\begin{align*}
\Da_{\bj}\nu^N(P_1\otimes\cdots\otimes P_\ell)=
\Da_{\bj}\I_{\frac{\hell}{2}}(P_1\otimes\cdots\otimes P_\ell)
+\frac{1}{N^2}\Da_{\bj}\I_{\frac{\hell}{2}+1}(P_1\otimes\cdots\otimes P_\ell)+\\
\cdots
+\frac{1}{N^{2g}}\Da_{\bj}\I_{\frac{\hell}{2}+g}(P_1\otimes\cdots\otimes P_\ell)
+o(\frac{1}{N^{2g}}).
\end{align*}
\end{prop}

Finally we prove Theorem \ref{theofree2},
\begin{theo}
For all
$\bj=(j_1,\cdots,j_n)$, for all $g\in \NN$, there exists $\eta>0$
such that for all $\bt$ in
$B_{\eta,c}$,
$$\Da_{\bj}F^N_{V_{\bt}}
=\Da_{\bj}F^0(\bt)+\cdots
+\Da_{\bj}F^g(\bt)+o(\frac{1}{N^{2g}}).$$
Besides, $\Da_{\bj}F^g$
is the generating function for rooted
maps of genus $g$ associated with $V$:
$$\Da_{\bj}F^g(\bt)=
\sum_{k_1,\cdots,k_n\in\NN}
\prod_i\frac{(-t_i)^{k_i}}{k_i!}
\C_g^{\bk+\bj}(q_{i_1},\cdots,q_{i_p})$$
where $\C_g^{k_1,\cdots,k_n}$ is the number of maps
on a surface
of
genus $g$ with
$k_i$ vertices of type $q_i$.
\end{theo}
\begin{dem}
The case $\bj=0$ is just Theorem \ref{theofree}.
Thus we can assume $\bj\neq 0$, for example $j_1\neq 0$.
Observe that for all $i$, 
$$\frac{\partial}{\partial t_i}F^N_{V_{\bt}}=-E[\mun(q_i)]=-\I(q_i)
-\frac{1}{N^2}\nu^N(q_i).$$
we can use the
Property \ref{propder}: there exists $\eta>0$ such that
for $\bt\in B_{\eta,c}$,
\begin{align*}
\Da_{\bj}F^N_{V_{\bt}}
&= -\Da_{\bj-{\mathbbm 1}_{i=1}}(\I(q_1)+\frac{1}{N^2}\nu^N(q_1))\\
&=-\Da_{\bj-{\mathbbm 1}_{i=1}}\I(q_1)-\frac{1}{N^2}
\Da_{\bj-{\mathbbm 1}_{i=1}}\nu^N(q_1)\\
&=-\Da_{\bj-{\mathbbm 1}_{i=1}}\I(q_1)-\frac{1}{N^2}
\Da_{\bj-{\mathbbm 1}_{i=1}}\I_1(q_1)-\cdots-\frac{1}{N^{2g}}
\Da_{\bj-{\mathbbm 1}_{i=1}}\I_g(q_1)\\
&+o(\frac{1}{N^{2g}})
\end{align*}
Observe now that according to Lemma \ref{dercomb}
\begin{align*}
\Da_{\bj-{\mathbbm 1}_{i=1}}\I_g(q_1)&=\sum_{k_1,\cdots,k_n\in\NN}
\prod_i\frac{(-t_i)^{k_i}}{k_i!}
\MM_g^{\bk+\bj-{\mathbbm 1}_{i=1}}(q_1)\\
&=-\sum_{k_1,\cdots,k_n\in\NN}
\prod_i\frac{(-t_i)^{k_i}}{k_i!}
\C_g^{\bk+\bj}
=-\Da_{\bj}F_g
\end{align*}
and by the same method,
$\Da_{\bj-{\mathbbm 1}_{i=1}}\I(q_1)=-\Da_{\bj}F_0$.
Thus we get,
$$\Da_{\bj}F^N_{V_{\bt}}=\Da_{\bj}F_0+\cdots+\frac{1}{N^{2g}}\Da_{\bj}F_g
+o(\frac{1}{N^{2g}}).$$
\end{dem}

\bibliography{ordresup}

\begin{thebibliography}{10}

\bibitem{ane}
{\sc An{\'e}, C., Blach{\`e}re, S., Chafa{\"{\i}}, D., Foug{\`e}res, P.,
  Gentil, I., Malrieu, F., Roberto, C., and Scheffer, G.}
\newblock {\em Sur les in\'egalit\'es de {S}obolev logarithmiques}, vol.~10 of
  {\em Panoramas et Synth\`eses [Panoramas and Syntheses]}.
\newblock Soci\'et\'e Math\'ematique de France, Paris, 2000.
\newblock With a preface by Dominique Bakry and Michel Ledoux.

\bibitem{BEH2}
{\sc Bertola, M., Eynard, B., and Harnad, J.}
\newblock Duality, biorthogonal polynomials and multi-matrix models.
\newblock {\em Comm. Math. Phys. 229}, 1 (2002), 73--120.

\bibitem{FGZ}
{\sc Di~Francesco P.~D., G.~P., and J., Z.-J.}
\newblock 2d gravity and random matrices.
\newblock {\em Phys. Rep.}, 254 (1995).

\bibitem{EML}
{\sc Ercolani, N.~M., and McLaughlin, K. D. T.-R.}
\newblock Asymptotics of the partition function for random matrices via
  {R}iemann-{H}ilbert techniques and applications to graphical enumeration.
\newblock {\em Int. Math. Res. Not.}, 14 (2003), 755--820.

\bibitem{EKK}
{\sc Eynard, B., Kokotov, A., and Korotkin, D.}
\newblock {$1/N\sp 2$}-correction to free energy in {H}ermitian two-matrix
  model.
\newblock {\em Lett. Math. Phys. 71}, 3 (2005), 199--207.

\bibitem{GCMP}
{\sc Guionnet, A.}
\newblock First order asymptotics of matrix integrals ; a rigorous approach
  towards the understanding of matrix models.
\newblock {\em Comm. Math. Phys.\/} (2003).

\bibitem{GMa}
{\sc Guionnet, A., and Maurel-Segala, E.}
\newblock Combinatorial aspects of matrix models.
\newblock {\em arviv front.math.ucdavis.edu/math.PR/0503064\/} (2005).

\bibitem{GMa2}
{\sc Guionnet, A., and Maurel-Segala, E.}
\newblock Second order asymptotics for matrix models.
\newblock {\em arviv front.math.ucdavis.edu/math.PR/0601040\/} (2006).

\bibitem{Tu}
{\sc Tutte, W.~T.}
\newblock On the enumeration of planar maps.
\newblock {\em Bull. Amer. Math. Soc. 74\/} (1968), 64--74.

\bibitem{zvon}
{\sc Zvonkin, A.}
\newblock Matrix integrals and map enumeration: an accessible introduction.
\newblock {\em Math. Comput. Modelling 26}, 8-10 (1997), 281--304.
\newblock Combinatorics and physics (Marseilles, 1995).

\end{thebibliography}
\bibliographystyle{acm}
\end{document}